\documentclass[reqno,12pt,a4paper]{amsart}

\usepackage{mathrsfs}
\usepackage{amssymb}
\usepackage{amsfonts}
\usepackage{latexsym}
\usepackage{amsthm}
\usepackage{graphicx}
\usepackage{xcolor}

\def\lb{\label}

\newcommand{\er}[1]{\textrm{(\ref{#1})}}

\begin{document}


\renewcommand{\theequation}{\arabic{section}.\arabic{equation}}
\theoremstyle{plain}
\newtheorem{theorem}{\bf Theorem}[section]
\newtheorem{lemma}[theorem]{\bf Lemma}
\newtheorem{corollary}[theorem]{\bf Corollary}
\newtheorem{proposition}[theorem]{\bf Proposition}
\newtheorem{definition}[theorem]{\bf Definition}

\theoremstyle{remark}
\newtheorem{remark}[theorem]{\bf Remark}
\newtheorem{example}[theorem]{\bf Example}

\def\a{\alpha}  \def\cA{{\mathcal A}}     \def\bA{{\bf A}}  \def\mA{{\mathscr A}}
\def\b{\beta}   \def\cB{{\mathcal B}}     \def\bB{{\bf B}}  \def\mB{{\mathscr B}}
\def\g{\gamma}  \def\cC{{\mathcal C}}     \def\bC{{\bf C}}  \def\mC{{\mathscr C}}
\def\G{\Gamma}  \def\cD{{\mathcal D}}     \def\bD{{\bf D}}  \def\mD{{\mathscr D}}
\def\d{\delta}  \def\cE{{\mathcal E}}     \def\bE{{\bf E}}  \def\mE{{\mathscr E}}
\def\D{\Delta}  \def\cF{{\mathcal F}}     \def\bF{{\bf F}}  \def\mF{{\mathscr F}}
\def\c{\chi}    \def\cG{{\mathcal G}}     \def\bG{{\bf G}}  \def\mG{{\mathscr G}}
\def\z{\zeta}   \def\cH{{\mathcal H}}     \def\bH{{\bf H}}  \def\mH{{\mathscr H}}
\def\e{\eta}    \def\cI{{\mathcal I}}     \def\bI{{\bf I}}  \def\mI{{\mathscr I}}
\def\p{\psi}    \def\cJ{{\mathcal J}}     \def\bJ{{\bf J}}  \def\mJ{{\mathscr J}}
\def\vT{\Theta} \def\cK{{\mathcal K}}     \def\bK{{\bf K}}  \def\mK{{\mathscr K}}
\def\k{\kappa}  \def\cL{{\mathcal L}}     \def\bL{{\bf L}}  \def\mL{{\mathscr L}}
\def\l{\lambda} \def\cM{{\mathcal M}}     \def\bM{{\bf M}}  \def\mM{{\mathscr M}}
\def\L{\Lambda} \def\cN{{\mathcal N}}     \def\bN{{\bf N}}  \def\mN{{\mathscr N}}
\def\m{\mu}     \def\cO{{\mathcal O}}     \def\bO{{\bf O}}  \def\mO{{\mathscr O}}
\def\n{\nu}     \def\cP{{\mathcal P}}     \def\bP{{\bf P}}  \def\mP{{\mathscr P}}
\def\r{\rho}    \def\cQ{{\mathcal Q}}     \def\bQ{{\bf Q}}  \def\mQ{{\mathscr Q}}
\def\s{\sigma}  \def\cR{{\mathcal R}}     \def\bR{{\bf R}}  \def\mR{{\mathscr R}}
\def\S{\Sigma}  \def\cS{{\mathcal S}}     \def\bS{{\bf S}}  \def\mS{{\mathscr S}}
\def\t{\tau}    \def\cT{{\mathcal T}}     \def\bT{{\bf T}}  \def\mT{{\mathscr T}}
\def\f{\phi}    \def\cU{{\mathcal U}}     \def\bU{{\bf U}}  \def\mU{{\mathscr U}}
\def\F{\Phi}    \def\cV{{\mathcal V}}     \def\bV{{\bf V}}  \def\mV{{\mathscr V}}
\def\P{\Psi}    \def\cW{{\mathcal W}}     \def\bW{{\bf W}}  \def\mW{{\mathscr W}}
\def\o{\omega}  \def\cX{{\mathcal X}}     \def\bX{{\bf X}}  \def\mX{{\mathscr X}}
\def\x{\xi}     \def\cY{{\mathcal Y}}     \def\bY{{\bf Y}}  \def\mY{{\mathscr Y}}
\def\X{\Xi}     \def\cZ{{\mathcal Z}}     \def\bZ{{\bf Z}}  \def\mZ{{\mathscr Z}}
\def\be{{\bf e}} \def\bc{{\bf c}}
\def\bv{{\bf v}} \def\bu{{\bf u}}
\def\Om{\Omega}
\def\bp{{\bf p}}\def\bq{{\bf q}}
\def\bx{{\bf x}} \def\by{{\bf y}}
\def\bbD{\pmb \Delta}
\def\mm{\mathrm m}
\def\mn{\mathrm n}

\newcommand{\mc}{\mathscr {c}}

\newcommand{\gA}{\mathfrak{A}}          \newcommand{\ga}{\mathfrak{a}}
\newcommand{\gB}{\mathfrak{B}}          \newcommand{\gb}{\mathfrak{b}}
\newcommand{\gC}{\mathfrak{C}}          \newcommand{\gc}{\mathfrak{c}}
\newcommand{\gD}{\mathfrak{D}}          \newcommand{\gd}{\mathfrak{d}}
\newcommand{\gE}{\mathfrak{E}}
\newcommand{\gF}{\mathfrak{F}}           \newcommand{\gf}{\mathfrak{f}}
\newcommand{\gG}{\mathfrak{G}}           
\newcommand{\gH}{\mathfrak{H}}           \newcommand{\gh}{\mathfrak{h}}
\newcommand{\gI}{\mathfrak{I}}           \newcommand{\gi}{\mathfrak{i}}
\newcommand{\gJ}{\mathfrak{J}}           \newcommand{\gj}{\mathfrak{j}}
\newcommand{\gK}{\mathfrak{K}}            \newcommand{\gk}{\mathfrak{k}}
\newcommand{\gL}{\mathfrak{L}}            \newcommand{\gl}{\mathfrak{l}}
\newcommand{\gM}{\mathfrak{M}}            \newcommand{\gm}{\mathfrak{m}}
\newcommand{\gN}{\mathfrak{N}}            \newcommand{\gn}{\mathfrak{n}}
\newcommand{\gO}{\mathfrak{O}}
\newcommand{\gP}{\mathfrak{P}}             \newcommand{\gp}{\mathfrak{p}}
\newcommand{\gQ}{\mathfrak{Q}}             \newcommand{\gq}{\mathfrak{q}}
\newcommand{\gR}{\mathfrak{R}}             \newcommand{\gr}{\mathfrak{r}}
\newcommand{\gS}{\mathfrak{S}}              \newcommand{\gs}{\mathfrak{s}}
\newcommand{\gT}{\mathfrak{T}}             \newcommand{\gt}{\mathfrak{t}}
\newcommand{\gU}{\mathfrak{U}}             \newcommand{\gu}{\mathfrak{u}}
\newcommand{\gV}{\mathfrak{V}}             \newcommand{\gv}{\mathfrak{v}}
\newcommand{\gW}{\mathfrak{W}}             \newcommand{\gw}{\mathfrak{w}}
\newcommand{\gX}{\mathfrak{X}}               \newcommand{\gx}{\mathfrak{x}}
\newcommand{\gY}{\mathfrak{Y}}              \newcommand{\gy}{\mathfrak{y}}
\newcommand{\gZ}{\mathfrak{Z}}             \newcommand{\gz}{\mathfrak{z}}

\def\ve{\varepsilon}   \def\vt{\vartheta}    \def\vp{\varphi}    \def\vk{\varkappa} \def\vr{\varrho}

\def\A{{\mathbb A}} \def\B{{\mathbb B}} \def\C{{\mathbb C}}
\def\dD{{\mathbb D}} \def\E{{\mathbb E}} \def\dF{{\mathbb F}} \def\dG{{\mathbb G}} \def\H{{\mathbb H}}\def\I{{\mathbb I}} \def\J{{\mathbb J}} \def\K{{\mathbb K}} \def\dL{{\mathbb L}}\def\M{{\mathbb M}} \def\N{{\mathbb N}} \def\O{{\mathbb O}} \def\dP{{\mathbb P}} \def\R{{\mathbb R}}\def\S{{\mathbb S}} \def\T{{\mathbb T}} \def\U{{\mathbb U}} \def\V{{\mathbb V}}\def\W{{\mathbb W}} \def\X{{\mathbb X}} \def\Y{{\mathbb Y}} \def\Z{{\mathbb Z}}


\def\la{\leftarrow}              \def\ra{\rightarrow}            \def\Ra{\Rightarrow}
\def\ua{\uparrow}                \def\da{\downarrow}
\def\lra{\leftrightarrow}        \def\Lra{\Leftrightarrow}


\def\lt{\biggl}                  \def\rt{\biggr}
\def\ol{\overline}               \def\wt{\widetilde}
\def\ul{\underline}
\def\no{\noindent}


\let\ge\geqslant                 \let\le\leqslant
\def\lan{\langle}                \def\ran{\rangle}
\def\/{\over}                    \def\iy{\infty}
\def\sm{\setminus}               \def\es{\emptyset}
\def\ss{\subset}                 \def\ts{\times}
\def\pa{\partial}                \def\os{\oplus}
\def\om{\ominus}                 \def\ev{\equiv}
\def\iint{\int\!\!\!\int}        \def\iintt{\mathop{\int\!\!\int\!\!\dots\!\!\int}\limits}
\def\el2{\ell^{\,2}}             \def\1{1\!\!1}
\def\sh{\sharp}
\def\wh{\widehat}
\def\bs{\backslash}
\def\intl{\int\limits}

\def\na{\mathop{\mathrm{\nabla}}\nolimits}
\def\sh{\mathop{\mathrm{sh}}\nolimits}
\def\ch{\mathop{\mathrm{ch}}\nolimits}
\def\where{\mathop{\mathrm{where}}\nolimits}
\def\all{\mathop{\mathrm{all}}\nolimits}
\def\as{\mathop{\mathrm{as}}\nolimits}
\def\Area{\mathop{\mathrm{Area}}\nolimits}
\def\arg{\mathop{\mathrm{arg}}\nolimits}
\def\const{\mathop{\mathrm{const}}\nolimits}
\def\det{\mathop{\mathrm{det}}\nolimits}
\def\diag{\mathop{\mathrm{diag}}\nolimits}
\def\diam{\mathop{\mathrm{diam}}\nolimits}
\def\dim{\mathop{\mathrm{dim}}\nolimits}
\def\dist{\mathop{\mathrm{dist}}\nolimits}
\def\Im{\mathop{\mathrm{Im}}\nolimits}
\def\Iso{\mathop{\mathrm{Iso}}\nolimits}
\def\Ker{\mathop{\mathrm{Ker}}\nolimits}
\def\Lip{\mathop{\mathrm{Lip}}\nolimits}
\def\rank{\mathop{\mathrm{rank}}\limits}
\def\Ran{\mathop{\mathrm{Ran}}\nolimits}
\def\Re{\mathop{\mathrm{Re}}\nolimits}
\def\Res{\mathop{\mathrm{Res}}\nolimits}
\def\res{\mathop{\mathrm{res}}\limits}
\def\sign{\mathop{\mathrm{sign}}\nolimits}
\def\span{\mathop{\mathrm{span}}\nolimits}
\def\supp{\mathop{\mathrm{supp}}\nolimits}
\def\Tr{\mathop{\mathrm{Tr}}\nolimits}
\def\BBox{\hspace{1mm}\vrule height6pt width5.5pt depth0pt \hspace{6pt}}


\newcommand\nh[2]{\widehat{#1}\vphantom{#1}^{(#2)}}
\def\dia{\diamond}

\def\Oplus{\bigoplus\nolimits}



\def\qqq{\qquad}
\def\qq{\quad}
\let\ge\geqslant
\let\le\leqslant
\let\geq\geqslant
\let\leq\leqslant
\newcommand{\ca}{\begin{cases}}
\newcommand{\ac}{\end{cases}}
\newcommand{\ma}{\begin{pmatrix}}
\newcommand{\am}{\end{pmatrix}}
\renewcommand{\[}{\begin{equation}}
\renewcommand{\]}{\end{equation}}
\def\eq{\begin{equation}}
\def\qe{\end{equation}}
\def\[{\begin{equation}}
\def\bu{\bullet}

\makeatletter
\@namedef{subjclassname@2020}{\textup{2020} Mathematics Subject Classification}
\makeatother

\title[Spectral invariants for Schr\"odinger operators on periodic graphs]{Spectral invariants for discrete Schr\"odinger operators on periodic graphs}

\date{\today}
\author[Natalia Saburova]{Natalia Saburova}
\address{Northern (Arctic) Federal University, Severnaya Dvina Emb. 17, Arkhangelsk, 163002, Russia,
 \ n.saburova@gmail.com, \ n.saburova@narfu.ru}

\subjclass[2020]{58J53, 35J10, 05C50}
\keywords{discrete periodic Schr\"odinger operators, periodic graphs, Floquet spectrum, spectral invariants, isospectral potentials}

\begin{abstract}
The aim of this article is to present a complete system of Floquet spectral invariants for the discrete Schr\"odinger operators with periodic potentials on periodic graphs. These invariants are polynomials in the potential and determined by cycles in the quotient graph from some specific cycle sets. We discuss some properties of these invariants and give an explicit expression for the linear and quadratic (in the potential) Floquet spectral invariants. The constructed system of spectral invariants can be used to study the sets of isospectral periodic potentials for the Schr\"odinger operators on periodic graphs. In particular, we deduce that under certain assumptions, if a real potential is isospectral to the zero (respectively, "degree") potential, then it must be the zero (respectively, "degree") potential.
\end{abstract}

\maketitle

\section {\lb{Sec1}Introduction}
\setcounter{equation}{0}
We consider discrete Schr\"odinger operators with periodic potentials on periodic graphs. Such operators are used in solid state physics to describe electron motion in a crystal in the tight binding approximation \cite{AM76}. The study of the spectra of these operators allows to explain many physical properties of crystals.

Using Floquet theory the spectral analysis of the discrete Schr\"odinger operator $H$ with a periodic potential on a $\G$-periodic graph $\cG$, where $\G$ is a lattice of rank $d$ in $\R^d$, can be reduced to studying spectra of a family of \emph{Floquet} operators $H(k)$ acting on the finite quotient graph $\cG_*=\cG/\G$ and depending on the parameter $k\in\T^d:=\R^d/(2\pi\Z)^d$ called \emph{quasimomentum}. For each $k\in\T^d$ the spectrum $\s(H(k))$ of the operator $H(k)$ consists of $\n$ eigenvalues (including multiplicities), where $\n$ is the number of the vertices of the quotient graph $\cG_*$. The family of the spectra $\s(H(k))$, $k\in\T^d$, is called the \emph{Floquet spectrum} of the Schr\"odinger operator $H$. The spectrum of $H(0)$ is called the \emph{periodic spectrum} of $H$.

The main goal of this paper is to present a complete system of Floquet (respectively, periodic) spectral invariants for the discrete Schr\"odinger operators with periodic potentials on arbitrary periodic graphs. This is a system of functionals depending on the potential and geometric parameters of the periodic graph (some characteristics of cycles in its quotient graph). The identity of these systems for two Schr\"odinger operators is equivalent to the identity of their Floquet (respectively, periodic) spectra. This system of spectral invariants can be used to study the sets of isospectral periodic potentials, i.e., periodic potentials for which the Schr\"odinger operator on a fixed periodic graph has the same Floquet (or periodic) spectrum. Spectral invariants for the discrete Schr\"odinger operators on the simplest periodic graph, the $d$-dimensional lattice $\Z^d$, $d\geq2$, and their applications to the problem of isospectral periodic potentials were studied by T.Kappeler in a series of papers \cite{Ka88a,Ka88b,Ka89} and by W.Liu \cite{L23,L24}. To the best of our knowledge there are no such kind of results for other periodic graphs (except for the triangular lattice also considered in \cite{L23}). We note that for the Schr\"odinger operator with a $\n$-periodic potential on the one-dimensional lattice $\Z$ the Floquet spectral invariants have no sense, due to the simple structure of the quotient graph, which is just a cycle on $\n$ vertices. For periodic graphs with a more complicated structure a family of the Floquet spectral invariants is quite rich (due to a rich set of cycles in their quotient graphs).

We mention that in the \emph{continuous case} a family of Floquet spectral invariants for the periodic Schr\"odinger operators in $\R^d$ was constructed in \cite{ERT84,MN86,V08}.

\medskip

The paper is organized as follows. Introduction contains definitions of periodic and quotient (fundamental) graphs. We also recall the notion of edge indices from \cite{KS14}. Then we define the discrete Schr\"odinger operator on periodic graphs and briefly describe its spectrum. In Section \ref{Sec2} we formulate our main results:

$\bu$ We construct a complete system of the Floquet spectral invariants $\cI_n^\mm(Q)$, $n=1,2,\ldots,\n$, $\mm\in\Z^d$, and a complete system of the periodic spectral invariants $\cI_n(Q)$, $n=1,2,\ldots,\n$, for the discrete Schr\"odinger operators with periodic potentials $Q$ on arbitrary periodic graphs (Theorem \ref{TSpI}). These invariants are obtained using the trace formulas for the Schr\"odinger operator from \cite{KS22}. The invariants $\cI_n^\mm(Q)$ are polynomials in the potential $Q$ and determined by cycles in the quotient graph of length at most $n-1$ and with \emph{index} $\mm$.

$\bu$ We express the periodic spectral invariants $\cI_n(Q)$ in terms of the Floquet spectral invariants $\cI_n^\mm(Q)$, $\mm\in\Z^d$ (Proposition \ref{spcN}). These expressions, in particular, yield that, in the contrast to the case of the one-dimensional lattice $\Z$, the periodic spectrum of the Schr\"odinger operator on a periodic graph, in general, does not determine its Floquet spectrum.

$\bu$ For a wide class of periodic graphs, we give explicit expressions for the first three periodic spectral invariants $\cI_n(Q)$, $n=1,2,3$, of the Schr\"odinger operator (Corollary \ref{CSpI}). As a simple consequence of these invariants, we deduce that if the fundamental graph $\cG_*$ has no loops, then a real potential isospectral to the zero (respectively, degree) potential must be the zero (respectively, degree) potential (Proposition \ref{Pzpo}).

$\bu$ We provide explicit expressions for the linear and quadratic (in the potential values) Floquet spectral invariants of the Schr\"odinger operator (Theorem \ref{TLQI}).

\no Finally, we illustrate the obtained results by some simple examples of periodic graphs. In particular, we consider the lattice $\Z^d$, $d\geq2$, and compare the obtained spectral invariants with the known ones from \cite{Ka89,L23,L24}.

In Section \ref{Sec3} we prove our results. First, we recall the trace formulas for the Schr\"odinger operator from \cite{KS22} (Theorem \ref{TPG}). Then, we prove Theorem \ref{TSpI} about complete systems of Floquet and periodic spectral invariants. Here we also present explicit expressions for the first periodic spectral invariants $\cI_n(Q)$ and Floquet spectral invariants $\cI_n^\mm(Q)$,  $n=1,2,3$, $\mm\in\Z^d$ (Proposition \ref{TrTO}) and prove Corollary \ref{CSpI}, Theorem \ref{TLQI}, and Proposition \ref{Pzpo}.

Section \ref{Sec4} is devoted to examples of the Schr\"odinger operators on some simple periodic graphs and their spectral invariants.

\subsection{Periodic graphs and edge indices.} \lb{ss1.1}
Let $\cG=(\cV,\cE)$ be a connected infinite graph embedded into the space $\R^d$, where $\cV$ is the set of its vertices and $\cE$ is the set of its unoriented edges. Considering each edge in $\cE$ to have two orientations, we introduce the set $\cA$ of all oriented edges. An edge starting at a vertex $u$ and ending at a vertex $v$ from $\cV$ will be denoted as the ordered pair $(u,v)\in\cA$. Let $\ol\be=(v,u)$ be the inverse edge of $\be=(u,v)\in\cA$. The vertices $u,v\in\cV$ will be called \emph{adjacent} and denoted by $u\sim v$, if $(u,v)\in\cA$. We define the degree ${\vk}_v$ of the vertex $v\in\cV$ as the number of all edges from $\cA$ starting at $v$.

Let $\G$ be a lattice of rank $d$ in $\R^d$ with a basis $\{\ga_1,\ldots,\ga_d\}$, i.e.,
$$
\G=\Big\{\ga\in\R^d: \ga=\sum_{s=1}^dn_s\ga_s, \; (n_s)_{s=1}^d\in\Z^d\Big\},
$$
and let
\[\lb{fuce}
\Omega=\Big\{\bx\in\R^d : \bx=\sum_{s=1}^dx_s\ga_s, \; (x_s)_{s=1}^d\in
[0,1)^d\Big\}
\]
be the \emph{fundamental cell} of the lattice $\G$. We define the equivalence relation on $\R^d$:
$$
\bx\equiv \by \; (\hspace{-4mm}\mod \G) \qq\Leftrightarrow\qq \bx-\by\in\G \qqq
\forall\, \bx,\by\in\R^d.
$$

We consider \emph{locally finite $\G$-periodic graphs} $\cG$, i.e., graphs satisfying the following conditions:
\begin{itemize}
  \item[1)] $\cG=\cG+\ga$ for any $\ga\in\G$, i.e., $\cG$ is invariant under translation by any vector $\ga\in\G$;
  \item[2)] the quotient graph  $\cG_*=\cG/\G$ is finite.
\end{itemize}
The basis vectors $\ga_1,\ldots,\ga_d$ of the lattice $\G$ are called the {\it periods}  of $\cG$. We also call the quotient graph $\cG_*=\cG/\G$ the \emph{fundamental graph} of the periodic graph $\cG$. The fundamental graph $\cG_*$ is a graph on the $d$-dimensional torus $\R^d/\G$. The graph $\cG_*=(\cV_*,\cE_*)$ has the vertex set $\cV_*=\cV/\G$, the set $\cE_*=\cE/\G$ of unoriented edges and the set $\cA_*=\cA/\G$ of doubled oriented edges.

\begin{remark}
A periodic graph $\cG$ could be defined as an abstract infinite graph  equipped with an action of a finitely generated free abelian group (a lattice) $\G$  (see, e.g., \cite[Chapter~4]{BK13}). The embedding of $\cG$ into $\R^d$ is just a simple geometric realization of periodic graphs.
\end{remark}

We define the important notion of an {\it edge index} which was introduced in \cite{KS14}. This notion allows one to consider, instead of a periodic graph, the finite fundamental graph with edges labeled by some integer vectors called indices.

For each $\bx\in\R^d$, we denote by $\bx_\A\in\R^d$ the coordinate vector of $\bx$ with respect to the basis  $\A=\{\ga_1,\ldots,\ga_d\}$ of the lattice~$\G$, i.e.,
\[\lb{cola}
\bx_\A=(x_1,\ldots,x_d), \qqq \textrm{where} \qq \bx=
\textstyle\sum\limits_{s=1}^dx_s\ga_s.
\]

For any vertex $v\in\cV$ of a $\G$-periodic graph $\cG$, the following  unique representation holds true:
\[\lb{Dv}
v=v_0+[v], \qqq \textrm{where}\qqq v_0\in\cV\cap\Omega,\qqq [v]\in\G,
\]
and $\Omega$ is the fundamental cell of the lattice $\G$ defined by \er{fuce}. In other words, each vertex $v$ can be obtained from a vertex $v_0$ in $\Omega$  by a shift by a vector $[v]$ of the lattice $\G$.

For any oriented edge $\be=(u,v)\in\cA$ of the periodic graph $\cG$, we define the \emph{edge index} $\t(\be)$ as the vector of the lattice $\Z^d$ given by
\[\lb{in}
\t(\be)=[v]_\A-[u]_\A\in\Z^d,
\]
where $[v]\in\G$ is defined by \er{Dv} and the vector $[v]_\A\in\Z^d$  is given by \er{cola}. Due to the periodicity of the graph $\cG$, the edge indices $\t(\be)$ satisfy
\[\lb{Gpe}
\t(\be+\ga)=\t(\be),\qqq \forall\, (\be,\ga)\in\cA \ts\G.
\]
This periodicity of the indices allows us to assign the \emph{index} $\t(\be_*)\in\Z^d$ to each oriented edge $\be_*\in\cA_*$ of the fundamental graph $\cG_*$ by setting
\[\lb{dco}
\t(\be_*)=\t(\be),
\]
where $\be\in\cA$ is an oriented edge in the periodic graph $\cG$ from the equivalence class $\be_*\in\cA_*=\cA/\G$. In other words, edge indices of the fundamental graph $\cG_*$ are induced by edge indices of the periodic graph~$\cG$. Due to (\ref{Gpe}), the edge index $\t(\be_*)$ is uniquely determined by \er{dco} and does not depend on the choice of $\be\in\cA$. From the definition of the edge indices, it follows that
\[\lb{inin}
\t(\ol\be\,)=-\t(\be), \qqq \forall\,\be\in\cA_*,
\]
where $\ol\be$ is the inverse edge of $\be\in\cA_*$.

\subsection{Schr\"odinger operators on periodic graphs}
Let $\ell^2(\cV)$ be the Hilbert space of all square summable
functions  $f:\cV\to \C$ equipped with the norm
$$
\|f\|^2_{\ell^2(\cV)}=\sum_{v\in\cV}|f(v)|^2<\infty.
$$

We consider the discrete Schr\"odinger operator $H$ acting on $\ell^2(\cV)$ and given by
\[
\lb{Sh} H=A+Q,
\]
where $A$ is the \emph{adjacency} operator having the form
\[\lb{ALO}
(A f)(v)=\sum_{(v,u)\in\cA}f(u), \qqq f\in\ell^2(\cV), \qqq v\in\cV,
\]
and $Q$ is a $\G$-periodic potential (possibly complex valued), i.e., it satisfies
\[\lb{ppot}
(Qf)(v)=Q(v)f(v), \qqq Q(v+\ga)=Q(v), \qqq \forall\,(v,\ga)\in\cV\ts\G.
\]
The sum in \er{ALO} is taken over all edges from $\cA$ starting at the
vertex $v$.

\begin{remark}\lb{Re12}
\emph{i}) It is more usual and physically motivated to consider the Schr\"odinger operator $H=\D+Q$, where $\D$ is the discrete \emph{combinatorial Laplacian} given by $\D=\vk-A$. Here $\vk$ is the degree potential, i.e.,
\[\lb{depo}
(\vk f)(v)=\vk_v f(v), \qqq f\in\ell^2(\cV), \qqq v\in\cV,
\]
where $\vk_v$ is the degree of the vertex $v$. But since we do not impose any restrictions on the potential $Q$ (except $\G$-periodicity), we may absorb the degree potential $\vk$ into $Q$. Moreover, if $\cG$ is a regular graph of degree $\vk_o$, i.e., all vertices of $\cG$ have the same degree $\vk_o$, then the Laplacian $\D$ has the form $\D=\vk_oI-A$, where $I$ is the identity operator. Thus, the operators $-\D$ and $A$ on a regular graph differ only by a shift.

\emph{ii}) When dealing with the combinatorial Laplacian, without loss of generality, we may assume that \textbf{there are no loops in the periodic graph} $\cG$, see \cite[Section 1.1 (Remark~2)]{KS22}.

\emph{iii}) If $\cG$ is a $\G$-periodic graph and $Q$ is a $\G'$-periodic  potential, where $\G'$ is a sublattice of $\G$, then instead of the fundamental graph $\cG_*=\cG/\G$ of $\cG$ we will deal with its \emph{expanded fundamental graph} $\cG_*=\cG/\G'$.
\end{remark}

\subsection{Spectra of Schr\"odinger operators}
We briefly describe the spectrum of the Schr\"odinger operator $H$ on periodic graphs (for more details, see, e.g., \cite{HN09} or \cite{KS14}). We introduce the Hilbert space
$$
\mH=L^2\Big(\T^{d},{dk\/(2\pi)^d}\,;\ell^2(\cV_*)\Big)
=\int_{\T^{d}}^{\os}\ell^2(\cV_*)\,{dk \/(2\pi)^d}\,, \qqq
\T^d=\R^d/(2\pi\Z)^d,
$$
equipped with the norm
$$
\|g\|^2_{\mH}=\int_{\T^d}\|g(k)\|_{\ell^2(\cV_*)}^2\frac{dk}{(2\pi)^d}\,,\qqq g\in\mH.
$$

The Schr\"odinger operator $H$ on $\ell^2(\cV)$ has the following decomposition into a constant fiber direct integral, see \cite[Theorem 1.1]{KS14},
$$
UH U^{-1}=\int^\oplus_{\T^d}H(k){dk\/(2\pi)^d}\,,
$$
where $U:\ell^2(\cV)\to\mH$ is some unitary operator (the Gelfand
transform). The parameter $k$ is called the \emph{quasimomentum}.
For each $k\in\T^d$, the \emph{fiber} (or \emph{Floquet}) operator $H(k)$ on $\ell^2(\cV_*)$ is given by
\[\label{Hvt'}
H(k)=A(k)+Q.
\]
Here $Q$ is the potential on $\ell^2(\cV_*)$, and $A(k)$ is the fiber adjacency operator having the form
\[
\label{fado}
\big(A(k)f\big)(v)=\sum_{\be=(v,u)\in\cA_*}e^{i\lan\t(\be),\,k\ran}f(u),
 \qqq f\in\ell^2(\cV_*),\qqq v\in \cV_*,
\]
where $\t(\be)\in\Z^d$ is the index of the edge $\be\in\cA_*$ defined by \er{in}, \er{dco}, and $\lan\cdot,\cdot\ran$ denotes the standard inner product in $\R^d$.

Let $\#M$ denote the number of elements in a set $M$. For each $k\in\T^d$, the spectrum $\s(H(k))$ of $H(k)$ consists of $\n:=\#\cV_*$ eigenvalues (including multiplicities). The family of spectra of all Floquet operators $H(k)$
$$
\big\{\s\big(H(k)\big): k\in\T^d\big\}
$$
is called the \emph{Floquet} or, in physics literature, \emph{Bloch} spectrum of the Schr\"odinger operator $H$. The Floquet operator $H(0)$ is just the Schr\"odinger operator defined by \er{Sh} -- \er{ppot} on the finite fundamental graph $\cG_*$. Its spectrum $\s\big(H(0)\big)$ is called the \emph{periodic} spectrum of $H$.

If the potential $Q$ is real, then each fiber operator $H(k)$, $k\in\T^{d}$, is self-adjoint and has $\n$ real eigenvalues $\l_j(k)$, $j=1,\ldots,\n$, labeled in non-decreasing order counting multiplicities:
$$
\l_{1}(k)\leq\l_{2}(k)\leq\ldots\leq\l_{\nu}(k), \qqq
\forall\,k\in\T^{d}.
$$
Each \emph{band function} $\l_j(\cdot)$ is a continuous
and piecewise real analytic function on the torus $\T^{d}$ and creates the
\emph{spectral band} $\s_j(H)$ given by
\[\lb{ban.1H}
\begin{array}{l}
\s_j(H)=[\l_j^-,\l_j^+]=\l_j(\T^{d}),\qqq j\in\N_\n=\{1,\ldots,\n\},\\[8pt]
\displaystyle\textrm{where}\qqq \l_j^-=\min_{k\in\T^d}\l_j(k),\qqq \l_j^+=\max_{k\in\T^d}\l_j(k).
\end{array}
\]
The spectrum of the Schr\"odinger operator $H$ with a real periodic potential $Q$ on a periodic graph $\cG$ has the form
$$
\s(H)=\bigcup_{k\in\T^d}\s\big(H(k)\big)=
\bigcup_{j=1}^{\nu}\s_j(H),
$$
i.e., it consists of $\n$ bands $\s_j(H)$ defined by \er{ban.1H}.

\begin{remark}
Some of spectral bands $\s_j(H)=[\l_j^-,\l_j^+]$ may be degenerate, i.e., $\l_j^-=\l_j^+$.
\end{remark}

\section{Main results}
\setcounter{equation}{0}
\lb{Sec2}
\subsection{Cycles and their characteristics} \lb{ssci}
The main role in the construction of spectral invariants for the Schr\"odinger operator $H$ on a periodic graph $\cG$ is played by cycles in the fundamental graph $\cG_*$ and some of their characteristics.

A \emph{path} $\gp$ in a graph $\cG=(\cV,\cA)$ is a sequence of oriented edges
$$
\gp=(\be_1,\be_2,\ldots,\be_n), \qqq \textrm{where} \qqq
\be_j=(v_j,v_{j+1})\in\cA,\qq j=1,\ldots,n,
$$
for some vertices $v_1,\ldots,v_n,v_{n+1}\in\cV$. The number $n$ of edges in a path $\gp$ is called the \emph{length} of $\gp$ and is denoted by $|\gp|$, i.e., $|\gp|=n$. A path $\gp$ is called \emph{closed}, if $v_1=v_{n+1}$.
Repeating a closed path $\gp$ $r$ times, we obtain \emph{$r$-multiple} $\gp^r$ of $\gp$. If $\gp$ is not an $r$-multiple of a closed path with $r\geq2$, then $\gp$ is called \emph{prime}.

Two closed paths are \emph{equivalent} if one is obtained by a cyclic permutation of edges in another. A (prime) \emph{cycle} is an equivalence class of a (prime) closed path $\gp=(\be_1,\be_2,\ldots,\be_n)$ and will be denoted by $\bc=[\gp]=[\be_1,\be_2,\ldots,\be_n]$, its length $|\bc|=n$. The \emph{reverse} of a cycle $\bc=[\be_1,\ldots,\be_n]$ is the cycle $\ol\bc=[\ol\be_n,\ldots,\ol\be_1]$, where $\ol\be$ is the inverse edge of $\be\in\cA$.

\begin{remark}\lb{Rcdv} If a graph $\cG$ has no multiple edges, then a cycle $\bc$ can equivalently be described by the sequence $\bc=[v_1,\ldots,v_n]$ of vertices it passes by.
\end{remark}

Let $\cC$ be the set of all cycles in the fundamental graph $\cG_*=(\cV_*,\cA_*)$.  For any cycle $\bc\in\cC$, we define the \emph{cycle index} $\t(\bc)\in\Z^d$ by
\[\lb{cyin}
\t(\bc)=\sum\limits_{j=1}^n\t(\be_j),  \qqq \bc=[\be_1,\be_2,\ldots,\be_n]\in\cC,
\]
where $\t(\be)\in\Z^d$ is the index of the edge $\be\in\cA_*$ defined by \er{in}, \er{dco}. From this definition and the identity \er{inin}, it follows that
\[\lb{ininc}
\t(\ol\bc\,)=-\t(\bc), \qqq \forall\,\bc\in\cC.
\]

\begin{remark}\lb{Rein}
\emph{i}) Edge indices depend on the choice of the embedding of the periodic graph $\cG$ into $\R^d$. Cycle indices \emph{do not} depend on this choice. Indeed, any cycle $\bc$ in the fundamental graph $\cG_*$ is obtained by factorization of a path in the periodic graph $\cG$ connecting some $\G$-equivalent vertices $v\in\cV$ and $v+\ga\in\cV$, $\ga\in\G$. The index $\t(\bc)$ of the cycle $\bc$ is just the coordinate vector of $\ga$ with respect to the basis $\{\ga_1,\ldots,\ga_d\}$ of the lattice $\G$, and, therefore, does not depend on the choice of the embedding. In particular, $\t(\bc)=0$ if and only if the cycle $\bc$ in $\cG_*$ corresponds to a cycle in~$\cG$.

\emph{ii}) The set $\cC$ includes cycles with \emph{back-tracking parts}, i.e., the cycles $[\be_1,\ldots,\be_n]$ for which $\be_{j+1}=\ol\be_j$ for some $j\in\N_n$ ($\be_{n+1}$ is understood as $\be_1$). In particular, each oriented edge $\be\in\cA_*$ produces a back-tracking cycle $\bc_{\be}=[\be,\ol\be\,]$ of length 2 (note that $\bc_{\be}=\bc_{\ol\be}$).
\end{remark}

The last cycle characteristic we need is the cycle weight. At each vertex $v$ of the fundamental graph $\cG_*=(\cV_*,\cA_*)$, we add a loop $\be_v$ with zero index and consider the \emph{modified} fundamental graph
\[\lb{mfg}
\wt\cG_*=(\cV_*,\wt\cA_*),
\qqq \wt\cA_*=\cA_*\cup\{\be_v\}_{v\in\cV_*},\qqq \t(\be_v)=0.
\]

Let $\wt\cC$ be the set of all cycles in $\wt\cG_*$. For each cycle $\bc\in\wt\cC$, we define the \emph{weight}
\[\lb{Wcy}
\begin{array}{l}
\o(\bc,Q)=\o(\be_1)\ldots\o(\be_n), \qqq \textrm{where}\qqq \bc=[\be_1,\ldots,\be_n]\in\wt\cC,\\[6pt]
\o(\be)=\left\{
\begin{array}{cl}
1,  & \qq \textrm{if} \qq  \be\in\cA_* \\[3pt]
Q(v), & \qq \textrm{if} \qq \be=\be_v
\end{array}\right..
\end{array}
\]
From this definition it follows that
\[\lb{inwe}
\o(\ol\bc,Q)=\o(\bc,Q), \qqq \forall\,\bc\in\wt\cC.
\]

\begin{remark}\lb{Rcyw}
\emph{i}) The weight $\o(\bc,Q)$ of a cycle $\bc\in\wt\cC$ is a monomial of degree at most $|\bc|$ in the potential values $Q(v)$ at the vertices $v$ of the cycle $\bc$.

\emph{ii}) For each cycle $\bc\in\cC$ of the fundamental graph $\cG_*$, the weight $\o(\bc,Q)=1$.
\end{remark}

\subsection{Spectral invariants for Schr\"odinger operators} Let $H=A+Q$ be the Schr\"odinger operator with a $\G$-periodic potential $Q$ on a $\G$-periodic graph $\cG$.

$\bu$  A functional $\cI(Q)$ is a \emph{Floquet spectral invariant} of $H=A+Q$, if $\cI(Q)$ has the following property: if the Floquet spectra of the operators $A+Q_1$ and $A+Q_2$ coincide, i.e.,
$$
\s\big(A(k)+Q_1\big)=\s\big(A(k)+Q_2\big) \qq \forall\, k\in\T^d,
$$
then $\cI(Q_1)=\cI(Q_2)$. Here $A(k)$ is the fiber adjacency operator defined by \er{fado}.

$\bu$  A functional $\cI(Q)$ is a \emph{periodic spectral invariant} of $H=A+Q$, if $\cI(Q)$ has the following property: if the periodic spectra of the operators $A+Q_1$ and $A+Q_2$ coincide, i.e.,
$$
\s\big(A(0)+Q_1\big)=\s\big(A(0)+Q_2\big),
$$
then $\cI(Q_1)=\cI(Q_2)$.

A system of functionals $\cI_s(Q)$, $s\in S$, is called a \emph{complete} system of Floquet (respectively, periodic) spectral invariants of the Schr\"odinger operator $H=A+Q$, if the equality $\cI_s(Q_1)=\cI_s(Q_2)$ for all values $s\in S$ implies that the Floquet (respectively, periodic) spectra of the operators $A+Q_1$ and $A+Q_2$ coincide, and conversely, i.e.,
$$
\s\big(A(k)+Q_1\big)=\s\big(A(k)+Q_2\big) \qq \forall\, k\in\T^d \qq \Leftrightarrow\qq  \cI_s(Q_1)=\cI_s(Q_2) \qq \forall\, s\in S
$$
$$
\textrm{(respectively,}\qq\s\big(A(0)+Q_1\big)=\s\big(A(0)+Q_2\big) \qq \Leftrightarrow\qq  \cI_s(Q_1)=\cI_s(Q_2) \qq \forall\, s\in S).
$$

\begin{remark}\lb{RTIn}
\emph{i}) Each periodic spectral invariant of the Schr\"odinger operator $H$ is also a Floquet spectral invariant of $H$, but the inverse is not true in general.

\emph{ii)} It is known (see, e.g., \cite[p.~44]{HJ85}) that $\n\ts\n$ matrices $A$ and $B$ have the same eigenvalues if and only if $\Tr A^n=\Tr B^n$ for all $n\in\N_\n$, where $\Tr A^n$ is the trace of $A^n$. Then the collections
$$
\big\{\Tr H^n(k),\; n\in\N_\n,\;k\in\T^d\big\}\qqq \textrm{and}\qqq \big\{\Tr H^n(0),\; n\in\N_\n\big\}
$$
are complete systems of Floquet and, respectively, periodic spectral invariants for the Schr\"o\-din\-ger operator $H$. Recall that $H(k)$ is the corresponding Floquet operator (matrix) defined by \er{Hvt'}, \er{fado}.
\end{remark}

In the following theorem we present complete systems of Floquet and periodic spectral invariants for the Schr\"odinger operator which are determined by cycles in the modified fundamental graph. Recall that each cycle $\bc$ has the following three characteristics:
\begin{itemize}
  \item $|\bc|$ is the length of the cycle $\bc$;
  \item $\t(\bc)$ is the index of $\bc$ defined by \er{cyin};
  \item $\o(\bc,Q)$ is the weight of $\bc$ defined by \er{Wcy}.
\end{itemize}

\begin{theorem}\lb{TSpI}
Let $H=A+Q$ be the Schr\"odinger operator on a periodic graph $\cG$. Then the following statements hold true.

i) The functionals
\[\lb{cInm}
\cI_n^\mm(Q)=\sum_{r\in\N_n,\;\bc\in\wt\cP\sm\cP\atop r|\bc|=n,\;r\t(\bc)=\mm}\frac1r\;\o^r(\bc,Q),\qqq n\in\N,\qqq \mm\in\Z^d,
\]
are Floquet spectral invariants of $H$. Here $\cP$ and $\wt\cP$ are the sets of all \textbf{prime} cycles in the fundamental graph $\cG_*=(\cV_*,\cA_*)$ and in the modified fundamental graph $\wt\cG_*$ defined by (\ref{mfg}), respectively. Moreover,
\[\lb{cssi}
\big\{\cI_n^\mm(Q),\; n\in\N_\n,\; \mm\in\Z^d\big\}, \qqq \nu=\#\cV_*,
\]
is a complete system of Floquet spectral invariants of $H$.

ii) The functionals
\[\lb{cIn}
\cI_n(Q)=\sum_{r\in\N_n,\;\bc\in\wt\cP\sm\cP\atop r|\bc|=n}\frac1r\;\o^r(\bc,Q),\qqq n\in\N,
\]
are periodic spectral invariants of $H$, and the system
\[\lb{css}
\big\{\cI_n(Q),\; n\in\N_\n\big\}
\]
is complete.
\end{theorem}

\begin{remark}
\emph{i}) For each $(n,\mm)\in\N\ts\Z^d$, the sum in \er{cInm} is taken over all $r\in\N_n=\{1,2,\ldots,n\}$ and all prime cycles $\bc\in\wt\cP\sm\cP$ such that $r|\bc|=n$ and $r\t(\bc)=\mm$. This condition yields that $r$ is a common divisor of $n$ and all components of $\mm$.

Similarly, for each $n\in\N$, the sum in \er{cIn} is taken over all $r\in\N_n$ and all prime cycles $\bc\in\wt\cP\sm\cP$ such that $r|\bc|=n$, i.e., $r$ is a divisor of $n$.

\emph{ii}) The sums in \er{cInm} and \er{cIn} are finite, since the number of cycles of length at most $n$ in the finite modified fundamental graph $\wt\cG_*$ is finite.

\emph{iii}) Each cycle $\bc\in\wt\cP\sm\cP$ has at least one loop edge $\be_v$, $v\in\cV_*$, added to the fundamental graph $\cG_*$, see \er{mfg}.
\end{remark}

We formulate some simple properties of the spectral invariants $\cI_n^\mm(Q)$ and $\cI_n(Q)$. A vector $\mm\in\Z^d$ with (setwise) coprime components is called \emph{primitive}.

\begin{proposition}\lb{spcN}
For the spectral invariants $\cI_n^\mm(Q)$ and $\cI_n(Q)$ defined by \er{cInm} and \er{cIn}, respectively, the following statements hold true.

i) $\cI_n^\mm(Q)$ and $\cI_n(Q)$ are polynomials in the potential values $Q(v)$, $v\in\cV_*$, of degree at most $n$, and they satisfy
\[\lb{coPF}
\cI_n(Q)=\sum_{\mm\in\Z^d}\cI_n^\mm(Q),\qqq n\in\N,
\]
\[\lb{Npme}
\cI_n^\mm(Q)=\cI_n^{-\mm}(Q),\qqq (n,\mm)\in\N\ts\Z^d,
\vspace{3mm}
\]
\[\lb{cNe0}
\cI_n^\mm(Q)=0, \qq \textrm{if}\qq
\|\mm\|>(n-1)\t_+,\qq\textrm{where}\qq\t_+=\max\limits_{\be\in\cA_*}\|\t(\be)\|,
\]
$\t(\be)$ is the index of the edge $\be\in\cA_*$ of the fundamental graph $\cG_*=(\cV_*,\cA_*)$ defined by \er{in}, \er{dco}, and $\|\cdot\|$ is the standard norm in $\R^d$.

ii) If $\mm\in\Z^d$ is primitive, then
\[\lb{Inpm}
\cI_n^\mm(Q)=\sum_{\bc\in\wt\cP_n^\mm\sm\cP_n^\mm}\o(\bc,Q),\qqq n\in\N,
\]
where $\cP_n^\mm$ and $\wt\cP_n^\mm$ are the sets of all prime cycles of length $n$ and with index $\mm$ in the fundamental graph $\cG_*$ and in the modified fundamental graph $\wt\cG_*$ defined by (\ref{mfg}), respectively; and $\o(\bc,Q)$ is the weight of a cycle $\bc$ defined by \er{Wcy}.

iii) If $n$ is prime, then
\[\lb{Ipnm}
\cI_n^\mm(Q)=\sum_{\bc\in\wt\cP_n^\mm\sm\cP_n^\mm}\o(\bc,Q),\qqq 0\neq\mm\in\Z^d,
\]
\[\lb{Ipn0}
\cI_n^0(Q)=\sum_{\bc\in\wt\cP_n^0\sm\cP_n^0}\o(\bc,Q)+\frac1n
\sum_{v\in\cV_*}Q^n(v),
\]
\[\lb{cIpn}
\cI_n(Q)=\sum_{\bc\in\wt\cP_n\sm\cP_n}\o(\bc,Q)+\frac1n
\sum_{v\in\cV_*}Q^n(v),
\]
where $\cP_n$ and $\wt\cP_n$ are the sets of all prime cycles of length $n$ in $\cG_*$ and $\wt\cG_*$, respectively.
\end{proposition}

\begin{remark}\lb{Re27}
\emph{i}) The number $\t_+$ in \er{cNe0} is the maximum "hopping" distance, i.e., the maximum $\Z^d$ distance between the copies of the fundamental domain whose vertices are connected by an edge.

\emph{ii}) The identity \er{cNe0} yields that the sum over $\mm$ in \er{coPF} is finite. From \er{cNe0} it also follows that the complete system of the Floquet spectral invariants \er{cssi} is finite (but not minimal). Due to \er{Npme}, the system of the functionals
$$
\cI_n^\mm(Q),\qq n\in\N_\n,\qq \mm=(m_1,\ldots,m_d)\in\Z^d,\qq m_1\geq0,
$$
is also a finite complete system of Floquet spectral invariants.

\emph{iii}) The proof of Theorem \ref{TSpI} is based on the trace formulas for the Schr\"odinger operator from \cite{KS22}.
\end{remark}

The next corollary gives simple explicit expressions for the first three periodic spectral invariants $\cI_n(Q)$ under certain conditions on the fundamental graph.

\begin{corollary}\lb{CSpI}
For $n=1,2,3$, the periodic spectral invariants $\cI_n(Q)$ defined by \er{cIn} have the form:\\[2pt]
$\bu$ $\cI_1(Q)=\sum\limits_{v\in\cV_*}Q(v)$,\\[2pt]
$\bu$ $\cI_2(Q)=\frac12\sum\limits_{v\in\cV_*}Q^2(v)$, if the fundamental graph $\cG_*=(\cV_*,\cA_*)$ has no loops,\\[2pt]
$\bu$ $\cI_3(Q)=\frac13\sum\limits_{v\in\cV_*}Q^3(v)+\vk_o\sum\limits_{v\in\cV_*}Q(v)$, if $\cG_*$ is a regular graph of degree $\vk_o$ without loops and multiple edges.
\end{corollary}

\begin{remark}
\emph{i}) For the  Schr\"odinger operator $H=A+Q$ with a $\G$-periodic potential $Q$ on the lattice $\Z^d$, where $\G=p_1\Z\oplus\ldots\oplus p_d\Z$ for integers $p_1,\ldots,p_d\geq2$, the periodic spectral invariants
\[\lb{ftpi}
\textstyle \sum\limits_{v\in\cV_*}Q(v),\qqq \sum\limits_{v\in\cV_*}Q^2(v),\qqq \sum\limits_{v\in\cV_*}Q^3(v)
\]
were derived in \cite{Ka89} using the discrete heat equation. Due to Corollary \ref{CSpI}, many periodic regular graphs (for example, the hexagonal lattice and the Kagome lattice with expanded fundamental graphs without multiple edges, see Fig.~\ref{fig2}\emph{c}) also have the periodic spectral invariants \er{ftpi}.

\emph{ii}) For (non-regular) fundamental graphs possibly having loops and multiple edges, the spectral invariants $\cI_2(Q)$ and $\cI_3(Q)$ are given in Proposition \ref{TrTO}.
\end{remark}

The following theorem provides the explicit expressions for the linear and quadratic (in the potential values) Floquet spectral invariants of the Schr\"odinger operator. A graph is called \emph{bipartite} if its vertex set is divided into two disjoint sets (called \emph{parts} of the graph) such that each edge connects vertices from distinct parts.

\begin{theorem}\lb{TLQI}
Let $\mm\in\Z^d$ be primitive, and $n:=n(\mm)$ be the length of the shortest cycle with index $\mm$ in the fundamental graph $\cG_*=(\cV_*,\cA_*)$. Then the Floquet spectral invariants $\cI^\mm_{n+s}(Q)$, $s=1,2$, defined by \er{cInm} have the form
\[\lb{lfin} \cI^\mm_{n+1}(Q)=\hspace{-3mm}\sum\limits_{\bc=[v_1,\ldots,v_n]\in\cP_n^\mm}
\hspace{-2mm}h_1(q_1,\ldots,q_n),
\]
\[\lb{qfin} \cI^\mm_{n+2}(Q)=\hspace{-3mm}\sum\limits_{\bc=[v_1,\ldots,v_n]\in\cP_n^\mm}
\hspace{-2mm}h_2(q_1,\ldots,q_n)+
\hspace{-3mm}\sum\limits_{\bc=[v_1,\ldots,v_{n+1}]\in\cP_{n+1}^\mm}
\hspace{-2mm}h_1(q_1,\ldots,q_{n+1}),
\]
where \\
$\bu$ $\bc=[v_1,\ldots,v_n]$ denotes a prime cycle $\bc$ of length $n$ given by the sequence of its vertices $v_1,\ldots,v_n\in\cV_*$; $q_j:=Q(v_j)$, $j\in\N_n$;\\
$\bu$ $\cP_n^\mm$ is the set of all prime cycles of length $n$ and with index $\mm$ in $\cG_*$; \\
$\bu$ $h_s(q_1,\ldots,q_n)$ is the complete homogeneous symmetric polynomial of degree $s$ in $q_1,\ldots,q_n$:
$$
h_1(q_1,\ldots,q_n)=\sum_{1\leq j\leq n} q_j,\qqq h_2(q_1,\ldots,q_n)=\sum_{1\leq j\leq l\leq n} q_jq_l\,.
$$

In particular, if the periodic graph $\cG$ is bipartite, then $\cP_{n+1}^\mm=\varnothing$ and
\[\lb{qbfi}
\cI^\mm_{n+2}(Q)=\hspace{-3mm}\sum\limits_{\bc=[v_1,\ldots,v_n]\in\cP_n^\mm}
\hspace{-2mm}h_2(q_1,\ldots,q_n).
\]
\end{theorem}

\begin{remark}\lb{Re2.12}
\emph{i}) In \er{lfin} -- \er{qbfi} the sum is taken over all \emph{distinct} prime cycles $\bc=[\be_1,\ldots,\be_n]\in\cP_n^\mm$, where $\be_j=(v_j,v_{j+1})$, $j\in\N_n$, $v_{n+1}=v_1$, even if some of them have the same sequence of vertices. This situation may happen if $\cG_*$ have multiple edges (see also Remark \ref{Rcdv}).

\emph{ii}) The Floquet spectral invariants \er{lfin} -- \er{qbfi}  are expressed in terms of only prime cycles of the fundamental graph $\cG_*$, not the modified fundamental graph $\wt\cG_*$ (compare with \er{cInm}).

\emph{iii}) Since $n(\mm)$ is the length of the shortest cycle with a primitive index $\mm$ in $\cG_*$, then $\cI^\mm_l(Q)=0$ for all $l\leq n(\mm)$.
\end{remark}

\subsection{Zero and degree potentials} Let $\cG$ be a $\G$-periodic graph. We say that $\G$-periodic potentials $Q_1$ and $Q_2$ are \emph{isospectral}, if the Schr\"odinger operators $H_1=A+Q_1$ and $H_2=A+Q_2$ on $\cG$ have the same \emph{periodic} spectrum.

The following proposition about isospectrality to the zero and degree potentials is an immediate consequence of the linear and quadratic periodic spectral invariants presented in Corollary \ref{CSpI}.

\begin{proposition}\lb{Pzpo}
Let $H=A+Q$ be the Schr\"odinger operator with a \textbf{real} $\G$-periodic potential $Q$ on a $\G$-periodic graph $\cG$.

i) Let $Q$ be isospectral to the zero potential, i.e., the Schr\"odinger operator $H=A+Q$ and the adjacency operator $A$ have the same periodic spectrum. If the fundamental graph $\cG_*=\cG/\G$ has no loops, then $Q\equiv0$.

ii) If $Q$ is isospectral to the minus degree potential $-\vk$ defined by \er{depo}, i.e., the Schr\"odinger operator $H=A+Q$ and the minus Laplacian $-\D:=A-\vk$ have the same periodic spectrum, then $Q\equiv-\vk$.
\end{proposition}

\begin{remark}

\emph{i}) For the lattice $\Z^d$ this result was proved in \cite{Ka89}. We suppose that for an arbitrary periodic graph it also should be known. However, we could not find it in the literature, which is the reason why we state it here.

\emph{ii}) In the proof of this proposition we follow arguments from \cite{ERT84} which were used to prove the analogous result in the continuous case.

\emph{iii}) If the fundamental graph $\cG_*$ has loops, then a non-zero \textit{real}  potential $Q$ isospectral to the zero potential may exist, see Example \ref{Enzp}.

\emph{iv}) Complex potentials $Q$ isospectral to the zero potential (or to the minus degree potential $-\vk$) exist. For an example of such complex potentials on $\Z^d$, see \cite{FLM24}.
\end{remark}

\subsection{Examples} In this section we obtain the spectral invariants $\cI_n^\mm(Q)$ and $\cI_n(Q)$ defined by \er{cInm} and \er{cIn} for the Schr\"odinger operator on some concrete periodic graphs. The proofs of the examples are based on Theorems \ref{TSpI}, \ref{TLQI}, and Corollary \ref{CSpI} and given in Section \ref{Sec4}.

\begin{figure}\centering
\unitlength 1.2mm
\begin{picture}(85,10)
\put(5,5.1){\line(1,0){75.00}} \put(10,5){\circle{1.2}}
\put(20,5){\circle*{1.2}} \put(30,5){\circle*{1.2}}
\put(40,5){\circle*{1.2}} \put(50,5){\circle*{1.2}}
\put(65,5){\circle*{1.2}} \put(75,5){\circle{1.2}}
\put(5,8){$\Z$} \put(9,1.0){$0$}
\put(19,1.0){$1$} \put(29,1.0){$2$} \put(39,1.0){$3$} \put(55,2.5){$\ldots$}
\put(49,1.0){$4$} \put(64,1.0){$\n$} \put(73,1.0){$\n+1$}
\put(2,1.0){\emph{a})}
\end{picture}

\unitlength 1.0mm
\begin{picture}(40,34)
\put(2,0){\emph{b})}
\put(4,23){$\cG_*$}
\put(6.5,13){$1$}
\put(15.5,24.5){$2$}
\put(29,20.1){$3$}
\put(29,6){$4$}
\put(21,2.5){$\ldots$}
\put(15.5,1.5){$\n$}
\put(10,14){\circle*{1.5}}
\put(17,23.2){\circle*{1.5}}
\put(17,4.8){\circle*{1.5}}
\put(27.6,8.1){\circle*{1.5}}
\put(27.6,19.9){\circle*{1.5}}
\bezier{600}(10,14)(11,23)(20,24)
\bezier{600}(30,14)(29,23)(20,24)
\bezier{600}(10,14)(11,5)(20,4)
\bezier{600}(30,14)(29,5)(20,4)
\end{picture}
\begin{picture}(40,34)
\put(2,0){\emph{c})}
\put(4,23){$\wt\cG_*$}
\put(6.5,13){$1$}
\put(15.5,24.5){$2$}
\put(29,20.1){$3$}
\put(29,6){$4$}
\put(21,2.5){$\ldots$}
\put(15.5,1.5){$\n$}
\put(10,14){\circle*{1.5}}
\put(17,23.2){\circle*{1.5}}
\put(17,4.8){\circle*{1.5}}
\put(27.6,8.1){\circle*{1.5}}
\put(27.6,19.9){\circle*{1.5}}
\bezier{600}(10,14)(11,23)(20,24)
\bezier{600}(30,14)(29,23)(20,24)
\bezier{600}(10,14)(11,5)(20,4)
\bezier{600}(30,14)(29,5)(20,4)
\color{red}
\put(5.8,14){\circle{8}}
\put(15.4,27.0){\circle{8}}
\put(15.4,0.8){\circle{8}}
\put(30.7,5.0){\circle{8}}
\put(30.7,22.9){\circle{8}}
\end{picture}
\caption{\scriptsize\emph{a}) The one-dimensional lattice $\Z$. \emph{b}) The fundamental graph $\cG_*=\Z/\n\Z$, $\n\in\N$. \emph{c})~The modified fundamental graph $\wt\cG_*$. The added loops are shown in red.}\label{fig3}
\end{figure}

\begin{example}\lb{E1DL}
Let $H=A+Q$ be the Schr\"odinger operator with a $\nu$-periodic potential $Q$ on the one-dimensional lattice $\Z$, $\n\in\N$, see Fig.~\ref{fig3}\emph{a}. Then the Floquet spectral invariants $\cI_n^\mm(Q)$ and the periodic spectral invariants $\cI_n(Q)$ of $H$ defined by \er{cInm} and \er{cIn}, respectively, satisfy
\[\lb{Inm1d}
\begin{array}{ll}
\cI_n^\mm(Q)=0, & \mm\in\Z\sm\{0\},\\[6pt]
\cI_n^0(Q)=\cI_n(Q),\qqq & n\in\N_\n.
\end{array}
\]
\end{example}

\begin{remark}
\emph{i}) For non-zero $\mm\in\Z$, the Floquet spectral invariants $\cI_n^\mm(Q)$, $n\in\N_\n$, have no sense. For each $n\in\N_\n$, the Floquet spectral invariant $\cI_n^0(Q)$ coincides with  the periodic spectral invariant $\cI_n(Q)$.

\emph{ii}) The Floquet spectrum of the Schr\"odinger operator $H$ on $\Z$ is determined by its periodic spectrum, see, e.g., \cite[p.~137]{Ka88a}. This fact also simply follows from \er{Inm1d} and completeness of the spectral invariants systems \er{cssi} and \er{css}.

\emph{iii}) The potential $Q$ on $\Z$ is not uniquely determined neither by the periodic spectrum, nor by the Floquet spectrum of $H$  \cite{Ka88a}. In particular, for a large \emph{generic} $\n$-periodic real potential $Q$ (i.e., a potential with pairwise distinct values at the fundamental graph vertices) there are $\n!$ real isospectral potentials.
\end{remark}

\begin{figure}[h]
\centering
\unitlength 1.0mm 
\linethickness{0.4pt}
\ifx\plotpoint\undefined\newsavebox{\plotpoint}\fi 
\begin{picture}(60,17)(0,0)
\put(5,5){\line(1,0){55.00}}
\put(10,5){\vector(1,0){15.00}}
\put(10,5){\line(0,1){10.00}}
\put(25,5){\line(0,1){10.00}}
\put(40,5){\line(0,1){10.00}}
\put(55,5){\line(0,1){10.00}}

\put(10,5){\circle*{1}}
\put(25,5){\circle{1}}
\put(40,5){\circle{1}}
\put(55,5){\circle{1}}

\put(10,15){\circle*{1}}
\put(25,15){\circle{1}}
\put(40,15){\circle{1}}
\put(55,15){\circle{1}}
\put(2,15){$\cG$}
\put(9,2){$v_1$}
\put(17,2){$\ga$}
\put(16,6){$\be_1$}
\put(24,2){$v_1\!\!+\!\ga$}
\put(9,16.5){$v_2$}
\put(22,16.5){$v_2\!+\!\ga$}
\put(6,10){$\be_2$}
\put(0,0){\emph{a})}
\end{picture}\hspace{20mm}
\begin{picture}(25,17)(0,0)
\put(10,5){\vector(0,1){10.00}}
\put(19.85,5.1){\vector(0,1){1.00}}
\put(15,5){\circle{10}}
\put(10,5){\circle*{1}}
\put(10,15){\circle*{1}}
\put(6,4.5){$v_1$}
\put(21.0,4){$\be_1$}
\put(6,10){$\be_2$}
\put(8.5,16.5){$v_2$}
\put(0,15){$\cG_*$}
\put(0,0){\emph{b})}
\end{picture}
\caption{\scriptsize \emph{a}) The graph $\cG$ is obtained by adding a pendant edge at each vertex of the one-dimensional lattice $\Z$; $\ga$ is the period of $\cG$\; \emph{b}) the fundamental graph $\cG_*=\cG/\Z$ consists of two vertices $v_1,v_2$ and two edges $\be_1,\be_2$ (and their inverse edges).} \label{slex1}
\end{figure}

\begin{example}\lb{Enzp}
Let $\cG$ be a graph obtained by adding a pendant edge at each vertex of the one-dimensional lattice $\Z$, see Fig.~\ref{slex1}\emph{a}. We consider the Schr\"odinger operator $H=A+Q$ with a $\Z$-periodic potential $Q$ on $\cG$. The fundamental graph $\cG_*=\cG/\Z$ has two vertices $v_1$ and $v_2$ with degrees $\vk_1=3$ and $\vk_2=1$, respectively (see Fig.~\ref{slex1}\emph{b}). Let $q_s:=Q(v_s)$, $s=1,2$. Then
\[\lb{1fsi}
\textstyle\cI_1^0(Q)=q_1+q_2,\qqq \cI_2^0(Q)=\frac12\,(q_1^2+q_2^2), \qqq \cI_2^1(Q)=\cI_2^{-1}(Q)=q_1
\]
form a complete system of Floquet spectral invariants of $H$; and
\[\lb{1psi}
\textstyle\cI_1(Q)=q_1+q_2,\qqq \cI_2(Q)=\frac12\,(q_1^2+q_2^2)+2q_1
\]
form a complete system of periodic spectral invariants of $H$.

If a (complex-valued) potential $P$ is isospectral to $Q$, then $P=Q=(q_1,q_2)$ or $P=(q_2-2,q_1+2)$. In particular, \\
$\bu$ if $Q=0$, then $P=0$ or $P=(-2,2)$; \\
$\bu$ if $Q=-\vk$, where $\vk=(\vk_1,\vk_2)=(3,1)$, then $P=-\vk$.
\end{example}

\begin{remark}
The previous example demonstrates that for the graph $\cG$ shown in Fig.~\ref{slex1}\emph{a}:

\emph{i}) The periodic spectrum does not determine the potential $Q$ uniquely. The Floquet spectrum does, see \er{1fsi}.

\emph{ii}) The periodic spectrum does not determine the Floquet spectrum as well. For example, the periodic spectrum of the Schr\"odinger operators $H=A+Q$ on the graph $\cG$ with $Q=0$ and $Q=(-2,2)$ is the same: $\{1\pm\sqrt{2}\,\}$, but their Floquet spectra are given by
$$
\cos k\pm\sqrt{\cos^2 k+1}\qqq\textrm{and}\qqq \cos k\pm\sqrt{\cos^2 k-4\cos k+5},\qqq k\in[-\pi,\pi],
$$
respectively.
\end{remark}

\begin{figure}[h]
\centering
\unitlength 1.1mm 
\linethickness{0.4pt}
\begin{picture}(80,58)(0,0)
\bezier{25}(21,15)(26,25)(31,35)
\bezier{25}(22,15)(27,25)(32,35)
\bezier{25}(23,15)(28,25)(33,35)
\bezier{25}(24,15)(29,25)(34,35)
\bezier{25}(25,15)(30,25)(35,35)
\bezier{25}(26,15)(31,25)(36,35)
\bezier{25}(27,15)(32,25)(37,35)
\bezier{25}(28,15)(33,25)(38,35)
\bezier{25}(29,15)(34,25)(39,35)
\bezier{25}(30,15)(35,25)(40,35)
\bezier{25}(31,15)(36,25)(41,35)
\bezier{25}(32,15)(37,25)(42,35)
\bezier{25}(33,15)(38,25)(43,35)
\bezier{25}(34,15)(39,25)(44,35)
\bezier{25}(35,15)(40,25)(45,35)
\bezier{25}(36,15)(41,25)(46,35)
\bezier{25}(37,15)(42,25)(47,35)
\bezier{25}(38,15)(43,25)(48,35)
\bezier{25}(39,15)(44,25)(49,35)

\put(8.0,51){$\cG$}
\put(-9.0,13){\emph{a})}
\put(-3,15){\line(1,0){66.0}}
\put(7,35){\line(1,0){66.0}}
\put(17,55){\line(1,0){66.0}}
\put(-1.5,12){\line(1,2){23.0}}
\put(18.5,12){\line(1,2){23.0}}
\put(38.5,12){\line(1,2){23.0}}
\put(58.5,12){\line(1,2){23.0}}

\put(11.5,12){\line(-1,2){8.0}}
\put(31.5,12){\line(-1,2){18.0}}
\put(51.5,12){\line(-1,2){23.0}}
\put(66.5,22){\line(-1,2){18.0}}
\put(76.5,42){\line(-1,2){8.0}}

\put(20,15){\vector(1,0){20.0}}
\put(20,15){\vector(1,2){10.0}}

\put(0,15){\circle{1.3}}
\put(10,15){\circle{1.3}}
\put(20,15){\circle*{1.3}}
\put(30,15){\circle*{1.3}}
\put(40,15){\circle{1.3}}
\put(50,15){\circle{1.3}}
\put(60,15){\circle{1.3}}

\put(5,25){\circle{1.3}}
\put(25,25){\circle*{1.3}}
\put(45,25){\circle{1.3}}
\put(65,25){\circle{1.3}}

\put(10,35){\circle{1.3}}
\put(20,35){\circle{1.3}}
\put(30,35){\circle{1.3}}
\put(40,35){\circle{1.3}}
\put(50,35){\circle{1.3}}
\put(60,35){\circle{1.3}}
\put(70,35){\circle{1.3}}

\put(15,45){\circle{1.3}}
\put(35,45){\circle{1.3}}
\put(55,45){\circle{1.3}}
\put(75,45){\circle{1.3}}

\put(20,55){\circle{1.3}}
\put(30,55){\circle{1.3}}
\put(40,55){\circle{1.3}}
\put(50,55){\circle{1.3}}
\put(60,55){\circle{1.3}}
\put(70,55){\circle{1.3}}
\put(80,55){\circle{1.3}}

\put(35.0,12.5){\small$\ga_1$}
\put(25.0,32.5){\small$\ga_2$}
\put(15.3,12.5){\small $v_1$}
\put(22.4,36.5){$\scriptstyle v_1+\ga_2$}
\put(39.5,36.0){$\scriptstyle v_2+\ga_2$}
\put(40.0,12.5){$\scriptstyle v_1+\ga_1$}
\put(20.5,25.5){\small$v_3$}
\put(26.3,45.0){$\scriptstyle v_3+\ga_2$}
\put(46.0,24.5){$\scriptstyle v_3+\ga_1$}
\put(29.5,16.5){\small$v_2$}
\put(33,24){$\Omega$}

\color{blue}
\put(24.0,16){\small$\be_1$}
\put(34.5,16){\small$\be_2$}
\put(20,15){\line(1,0){20.0}}
\color{red}
\put(18.5,20){\small$\be_3$}
\put(28.0,29){\small$\be_4$}
\put(20,15){\line(1,2){10.0}}
\color{brown}
\put(35,45){\line(1,-2){10.0}}
\put(38,40){\small$\be_5$}
\put(38.5,29){\small$\be_6$}
\end{picture}

\unitlength 1.0mm
\begin{picture}(40,50)(0,0)
\put(-5.0,10){\emph{b})}
\put(15.0,54){$\cG_*$}
\put(15,30){\circle*{1.5}}
\put(10.0,28.5){$v_1$}
\put(41,30){\circle*{1.5}}
\put(42,28.5){$v_2$}
\put(28,56){\circle*{1.5}}
\put(28.0,57){$v_3$}
\put(3,20){\footnotesize$\color{blue}\bc_1=[\be_1,\be_2]$, \qq $\t({\color{blue}\bc_1})=(1,0)$}
\put(3,15){\footnotesize$\color{red}\bc_2=[\be_3,\be_4]$, \qq $\t({\color{red}\bc_2})=(0,1)$}
\put(3,10){\footnotesize$\color{brown}\bc_3=[\be_5,\be_6]$, \qq $\t({\color{brown}\bc_3})=(1,-1)$}
\color{blue}
\bezier{200}(15,30)(28,18)(41,30)
\put(15,30){\line(1,0){26.0}}
\put(26.5,31.5){$\be_1$}
\put(26.5,25.5){$\be_2$}
\put(28,30.0){\vector(1,0){1.0}}
\put(28,24.0){\vector(-1,0){1.0}}

\color{red}
\put(16.3,46.5){\vector(-1,-2){1.0}}
\put(21.2,42.2){\vector(1,2){1.0}}
\bezier{200}(15,30)(10,48)(28,56)
\put(15,30){\line(1,2){13.0}}
\put(22.2,41){$\be_3$}
\put(10.0,45){$\be_4$}
\color{brown}
\put(41,30){\line(-1,2){13.0}}
\put(29.8,41){$\be_5$}
\put(41.5,45){$\be_6$}
\bezier{200}(41,30)(46,48)(28,56)
\put(34.2,43.5){\vector(1,-2){1.0}}
\put(40.6,44.5){\vector(-1,2){1.0}}
\end{picture}\hspace{10mm}
\begin{picture}(60,50)(0,0)
\put(18.0,54){$\cG_*$}
\put(13.0,10){\emph{c})}
\put(35,15){\circle*{1.3}}
\put(20,25){\circle*{1.3}}
\put(30,25){\circle*{1.3}}
\put(40,25){\circle*{1.3}}
\put(50,25){\circle*{1.3}}
\put(25,35){\circle*{1.3}}
\put(45,35){\circle*{1.3}}
\put(20,45){\circle*{1.3}}
\put(30,45){\circle*{1.3}}
\put(40,45){\circle*{1.3}}
\put(50,45){\circle*{1.3}}
\put(35,55){\circle*{1.3}}

\bezier{200}(20,25)(35,-1)(50,25)
\bezier{200}(20,45)(35,71)(50,45)

\bezier{200}(20,25)(6,52)(35,55)
\bezier{200}(50,25)(64,52)(35,55)
\bezier{200}(20,45)(6,18)(35,15)
\bezier{200}(50,45)(64,18)(35,15)

\put(20,25){\line(1,0){30.0}}
\put(20,45){\line(1,0){30.0}}
\put(20,25){\line(1,2){15.0}}
\put(35,15){\line(1,2){15.0}}
\put(35,15){\line(-1,2){15.0}}
\put(50,25){\line(-1,2){15.0}}

\end{picture}

\vspace{-10mm}
\caption{\scriptsize\emph{a}) The Kagome lattice $\cG$; $\ga_1,\ga_2$ are the periods of $\cG$. \emph{b}) The fundamental graph $\cG_*=\cG/\G$ of the Kagome lattice, where $\G$ is the lattice with the basis $\ga_1,\ga_2$; $\t(\bc_s)$ is the index of the cycle $\bc_s$, $s=1,2,3$. \emph{c}) The expanded fundamental graph $\cG_*=\cG/\G'$ of $\cG$, where $\G'$ is the lattice with the basis $2\ga_1,2\ga_2$.}\label{fig2}
\end{figure}

\begin{example}\lb{EKaL}
Let $H=A+Q$ be the Schr\"odinger operator with a $\G$-periodic potential $Q$ on the Kagome lattice $\cG$, where $\G$ is the lattice with the basis $\ga_1,\ga_2$, see Fig.~\ref{fig2}\emph{a}. Then
\[\lb{piKl}
\textstyle\cI_1(Q)=\sum\limits_{s=1}^3q_s,\qqq \cI_2(Q)=\frac12\sum\limits_{s=1}^3q^2_s,\qqq \cI_3(Q)=\frac13\sum\limits_{s=1}^3q^3_s+8\sum\limits_{s=1}^3q_s
\]
form a complete system of periodic spectral invariants of $H$; and
\[\lb{fiKl}
\cI_3^{(1,0)}=q_1+q_2,\qqq \cI_3^{(0,1)}=q_1+q_3,\qqq \cI_3^{(1,-1)}=q_2+q_3
\]
form a complete system of Floquet spectral invariants of $H$. Here $q_s:=Q(v_s)$, $s=1,2,3$, and $v_1,v_2,v_3$ are the vertices of the fundamental graph $\cG_*$ of the Kagome lattice (Fig.~\ref{fig2}\emph{b}).
\end{example}

\begin{remark}
\emph{i}) From \er{piKl} it follows that
\[\lb{piKl1}
\textstyle\sum\limits_{s=1}^3q_s,\qqq \sum\limits_{s=1}^3q^2_s,\qqq \sum\limits_{s=1}^3q^3_s
\]
also form a complete system of periodic spectral invariants of $H$.

\emph{ii}) The complete system of the Floquet spectral invariants (\ref{fiKl}) determines the periodic potential $Q$ uniquely. We note that the complete system of the periodic spectral invariants (\ref{piKl1}) gives (generically) $3!=6$ potentials with the same periodic spectrum.
\end{remark}

\begin{figure}[h]
\centering
\unitlength 1.0mm 
\linethickness{0.4pt}
\begin{picture}(45,47)(0,0)
\put(3,5){\line(1,0){44.00}}
\put(3,15){\line(1,0){44.00}}
\put(3,25){\line(1,0){44.00}}
\put(3,35){\line(1,0){44.00}}
\put(3,45){\line(1,0){44.00}}
\put(5,3){\line(0,1){44.00}}
\put(15,3){\line(0,1){44.00}}
\put(25,3){\line(0,1){44.00}}
\put(35,3){\line(0,1){44.00}}
\put(45,3){\line(0,1){44.00}}

\bezier{50}(5.5,5)(5.5,20)(5.5,35)
\bezier{50}(6.0,5)(6.0,20)(6.0,35)
\bezier{50}(6.5,5)(6.5,20)(6.5,35)
\bezier{50}(7.0,5)(7.0,20)(7.0,35)
\bezier{50}(7.5,5)(7.5,20)(7.5,35)
\bezier{50}(8.0,5)(8.0,20)(8.0,35)
\bezier{50}(8.5,5)(8.5,20)(8.5,35)
\bezier{50}(9.0,5)(9.0,20)(9.0,35)
\bezier{50}(9.5,5)(9.5,20)(9.5,35)
\bezier{50}(10.0,5)(10.0,20)(10.0,35)
\bezier{50}(10.5,5)(10.5,20)(10.5,35)
\bezier{50}(11.0,5)(11.0,20)(11.0,35)
\bezier{50}(11.5,5)(11.5,20)(11.5,35)
\bezier{50}(12.0,5)(12.0,20)(12.0,35)
\bezier{50}(12.5,5)(12.5,20)(12.5,35)
\bezier{50}(13.0,5)(13.0,20)(13.0,35)
\bezier{50}(13.5,5)(13.5,20)(13.5,35)
\bezier{50}(14.0,5)(14.0,20)(14.0,35)
\bezier{50}(14.5,5)(14.5,20)(14.5,35)
\bezier{50}(15.0,5)(15.0,20)(15.0,35)
\bezier{50}(15.5,5)(15.5,20)(15.5,35)
\bezier{50}(16.0,5)(16.0,20)(16.0,35)
\bezier{50}(16.5,5)(16.5,20)(16.5,35)
\bezier{50}(17.0,5)(17.0,20)(17.0,35)
\bezier{50}(17.5,5)(17.5,20)(17.5,35)
\bezier{50}(18.0,5)(18.0,20)(18.0,35)
\bezier{50}(18.5,5)(18.5,20)(18.5,35)
\bezier{50}(19.0,5)(19.0,20)(19.0,35)
\bezier{50}(19.5,5)(19.5,20)(19.5,35)
\bezier{50}(20.0,5)(20.0,20)(20.0,35)
\bezier{50}(20.5,5)(20.5,20)(20.5,35)
\bezier{50}(21.0,5)(21.0,20)(21.0,35)
\bezier{50}(21.5,5)(21.5,20)(21.5,35)
\bezier{50}(22.0,5)(22.0,20)(22.0,35)
\bezier{50}(22.5,5)(22.5,20)(22.5,35)
\bezier{50}(23.0,5)(23.0,20)(23.0,35)
\bezier{50}(23.5,5)(23.5,20)(23.5,35)
\bezier{50}(24.0,5)(24.0,20)(24.0,35)
\bezier{50}(24.5,5)(24.5,20)(24.5,35)
\bezier{50}(25.0,5)(25.0,20)(25.0,35)
\bezier{50}(25.5,5)(25.5,20)(25.5,35)
\bezier{50}(26.0,5)(26.0,20)(26.0,35)
\bezier{50}(26.5,5)(26.5,20)(26.5,35)
\bezier{50}(27.0,5)(27.0,20)(27.0,35)
\bezier{50}(27.5,5)(27.5,20)(27.5,35)
\bezier{50}(28.0,5)(28.0,20)(28.0,35)
\bezier{50}(28.5,5)(28.5,20)(28.5,35)
\bezier{50}(29.0,5)(29.0,20)(29.0,35)
\bezier{50}(29.5,5)(29.5,20)(29.5,35)
\bezier{50}(30.0,5)(30.0,20)(30.0,35)
\bezier{50}(30.5,5)(30.5,20)(30.5,35)
\bezier{50}(31.0,5)(31.0,20)(31.0,35)
\bezier{50}(31.5,5)(31.5,20)(31.5,35)
\bezier{50}(32.0,5)(32.0,20)(32.0,35)
\bezier{50}(32.5,5)(32.5,20)(32.5,35)
\bezier{50}(33.0,5)(33.0,20)(33.0,35)
\bezier{50}(33.5,5)(33.5,20)(33.5,35)
\bezier{50}(34.0,5)(34.0,20)(34.0,35)
\bezier{50}(34.5,5)(34.5,20)(34.5,35)

\put(5,5){\circle*{1.3}}
\put(15,5){\circle*{1.3}}
\put(25,5){\circle*{1.3}}
\put(35,5){\circle{1.3}}
\put(45,5){\circle{1.3}}
\put(5,15){\circle*{1.3}}
\put(15,15){\circle*{1.3}}
\put(25,15){\circle*{1.3}}
\put(35,15){\circle{1.3}}
\put(45,15){\circle{1.3}}
\put(5,25){\circle*{1.3}}
\put(15,25){\circle*{1.3}}
\put(25,25){\circle*{1.3}}
\put(35,25){\circle{1.3}}
\put(45,25){\circle{1.3}}
\put(5,35){\circle{1.3}}
\put(15,35){\circle{1.3}}
\put(25,35){\circle{1.3}}
\put(35,35){\circle{1.3}}
\put(45,35){\circle{1.3}}
\put(5,45){\circle{1.3}}
\put(15,45){\circle{1.3}}
\put(25,45){\circle{1.3}}
\put(35,45){\circle{1.3}}
\put(45,45){\circle{1.3}}
\put(5,5){\vector(1,0){30.00}}
\put(5,5){\vector(0,1){30.00}}
\put(5.5,6){$\scriptstyle(0,0)$}
\put(15.5,6){$\scriptstyle(1,0)$}
\put(25.5,6){$\scriptstyle(2,0)$}
\put(5.5,16){$\scriptstyle(0,1)$}
\put(15.5,16){$\scriptstyle(1,1)$}
\put(25.5,16){$\scriptstyle(2,1)$}
\put(5.5,26){$\scriptstyle(0,2)$}
\put(15.5,26){$\scriptstyle(1,2)$}
\put(25.5,26){$\scriptstyle(2,2)$}

\put(-1,39){$\Z^2$}
\put(31,2.5){$\ga_1$}
\put(1,32){$\ga_2$}
\put(19,19){$\Omega$}
\put(-3,3){\emph{a})}
\end{picture}\hspace{15mm} \unitlength 1.1mm
\begin{picture}(30,30)(0,0)
\put(-5,3){\emph{b})}
\put(5,15){\circle*{1.3}}
\put(15,15){\circle*{1.3}}
\put(25,15){\circle*{1.3}}
\put(5,25){\circle*{1.3}}
\put(15,25){\circle*{1.3}}
\put(25,25){\circle*{1.3}}
\put(5,5){\circle*{1.3}}
\put(15,5){\circle*{1.3}}
\put(25,5){\circle*{1.3}}
\put(2.1,2.7){$\scriptstyle(0,0)$}
\put(15.5,2.7){$\scriptstyle(1,0)$}
\put(25.5,2.7){$\scriptstyle(2,0)$}
\put(0.6,16.1){$\scriptstyle(0,1)$}
\put(15.5,16){$\scriptstyle(1,1)$}
\put(25.0,16){$\scriptstyle(2,1)$}
\put(0.2,26){$\scriptstyle(0,2)$}
\put(15.5,26){$\scriptstyle(1,2)$}
\put(25.5,26){$\scriptstyle(2,2)$}
\put(12,35){$\cG_*$}

\color{red}
\put(15,0.5){\vector(-1,0){1.00}}
\put(15,10.5){\vector(-1,0){1.00}}
\put(15,29.5){\vector(-1,0){1.00}}
\put(5,5){\line(1,0){20.00}}
\put(5,15){\line(1,0){20.00}}
\put(5,25){\line(1,0){20.00}}

\bezier{200}(5,15)(15,6)(25,15)
\bezier{200}(5,25)(15,34)(25,25)

\bezier{200}(5,5)(15,-4)(25,5)

\color{blue}
\put(10.5,15){\vector(0,-1){1.00}}
\put(29.5,15){\vector(0,-1){1.00}}
\put(0.5,15){\vector(0,-1){1.00}}
\bezier{200}(15,5)(6,15)(15,25)
\bezier{200}(25,5)(34,15)(25,25)
\bezier{200}(5,5)(-4,15)(5,25)
\put(5,5){\line(0,1){20.00}}
\put(15,5){\line(0,1){20.00}}
\put(25,5){\line(0,1){20.00}}
\end{picture}\hspace{18mm} \unitlength 0.8mm
\begin{picture}(30,30)(0,0)
\put(-9,3){\emph{c})}
\put(5,25){\circle*{1.7}}
\put(25,25){\circle*{1.7}}
\put(5,5){\circle*{1.7}}
\put(25,5){\circle*{1.7}}
\put(-3,2.0){$\scriptstyle(0,0)$}
\put(25.5,2.0){$\scriptstyle(1,0)$}
\put(-3,26){$\scriptstyle(0,1)$}
\put(26.0,26){$\scriptstyle(1,1)$}
\put(12,35){$\cG_*$}

\color{red}
\put(15,0.5){\vector(-1,0){1.00}}
\put(15,29.5){\vector(-1,0){1.00}}
\put(5,5){\line(1,0){20.00}}
\put(5,25){\line(1,0){20.00}}
\bezier{200}(5,25)(15,34)(25,25)
\bezier{200}(5,5)(15,-4)(25,5)

\color{blue}
\put(0.5,15){\vector(0,-1){1.00}}
\put(29.5,15){\vector(0,-1){1.00}}
\bezier{200}(25,5)(34,15)(25,25)
\bezier{200}(5,5)(-4,15)(5,25)
\put(5,5){\line(0,1){20.00}}
\put(25,5){\line(0,1){20.00}}
\end{picture}
\caption{\scriptsize \emph{a}) The square lattice $\Z^2$. \emph{b}) The fundamental graph $\cG_*=\Z^2/\G$ of the square lattice, where $\G=3\Z\oplus3\Z$. The vectors $\ga_1=(3,0)$, $\ga_2=(0,3)$ are the basis of the lattice $\G$, the fundamental cell $\Omega$ of $\G$ is shaded. The red cycles in $\cG_*$ have the index $(1,0)$ and the index of the blue cycles is $(0,1)$. \emph{c}) The fundamental graph $\cG_*=\Z^2/\G$ of the square lattice, where $\G=2\Z\oplus2\Z$.}\label{fig1}
\end{figure}

\begin{example}\lb{EsqL}
Let $H=A+Q$ be the Schr\"odinger operator with a $\G$-periodic potential $Q$ on the lattice $\Z^d$ (for $d=2$ see Fig.~\ref{fig1}\emph{a}), where $\G=p_1\Z\oplus\ldots\oplus p_d\Z$ for some integers $p_1,\ldots,p_d\geq2$. Then

\emph{i}) the periodic spectral invariants $\cI_n(Q)$, $n=1,2,3$, defined by \er{cIn} have the form
\[\lb{pisl1}
\cI_1(Q)=\sum_{\mn\in\cV_*}Q(\mn),\qqq \cI_2(Q)=\frac12\sum_{\mn\in\cV_*}Q^2(\mn),
\]
\[\lb{pisl2}
\cI_3(Q)=\frac13\sum\limits_{\mn\in\cV_*}Q^3(\mn)+2\,(d+d_o)\,\cI_1(Q),
\]
where $\cV_*$ is the vertex set of the fundamental graph $\cG_*=\Z^d/\G$ of the lattice $\Z^d$:
\[\lb{cV*}
\cV_*=\{\mn=(n_1,\ldots,n_d)\in\Z^d : 0\leq n_j<p_j, \; j=1,\ldots,d\},
\]
and $d_o$ is the number of integers $p_1,\ldots,p_d$ equal to 2;

\emph{ii}) the linear and quadratic Floquet spectral invariants $\cI^\mm_{n+s}(Q)$, $s=1,2$, given by \er{lfin} and \er{qbfi}, for $\mm=\mathfrak{e}_j$, and $n(\mm)=p_j$, $j=1,\ldots,d$, where $\{\mathfrak{e}_j\}_{j=1}^d$ is the standard basis of $\Z^d$, have the form
\[\lb{lfZd}
\cI^{\mathfrak{e}_j}_{p_j+1}(Q)=\cI_1(Q),
\]
\[\lb{qfZd}
\cI^{\mathfrak{e}_j}_{p_j+2}(Q)=\cI_2(Q)+\frac12\sum\limits_{n_1,\ldots,\breve{n}_j,\ldots,n_d}
\Big(\sum\limits_{n_j=0}^{p_j-1}Q(\mn)\Big)^2, \qq \mn=(n_1,\ldots,n_d),
\]
where $\cI_1(Q)$ and $\cI_2(Q)$ are given in \er{pisl1}, and the sum $\sum\limits_{n_1,\ldots,\breve{n}_j,\ldots,n_d}$ in \er{qfZd} is taken over all $n_1,\ldots,n_d$, except $n_j$, from 0 to ${p_1-1},\ldots,{p_d-1}$, respectively.
\end{example}

\begin{remark}
\emph{i}) In \er{cV*} we identify the vertices of the fundamental graph $\cG_*=\Z^d/\G$ with the corresponding vertices of $\Z^d$ in the fundamental cell $\Omega$ of the period lattice $\G$, see Fig.~\ref{fig1}.

\emph{ii}) The linear Floquet spectral invariants \er{lfZd} coincide with the periodic spectral invariant $\cI_1(Q)$ given in \er{pisl1}.

\emph{iii}) For each $j$, the sum in the round brackets in \er{qfZd} is the \emph{potential over a cycle} (i.e., the sum of the potential values at the vertices of the cycle) with index $\mm=\mathfrak{e}_j$ in the fundamental graph $\cG_*=\Z^d/\G$. For example, if $d=2$ and $p_1=p_2=3$ (see Fig.~\ref{fig1}\emph{b}), the Floquet spectral invariants \er{qfZd} have the form
\[\lb{fisl3}
\begin{array}{l}
\cI^{(1,0)}_5(Q)=\cI_2(Q)+\frac12\sum\limits_{n_2=0}^2
\Big(\sum\limits_{n_1=0}^2Q(n_1,n_2)\Big)^2,\\
\cI^{(0,1)}_5(Q)=\cI_2(Q)+\frac12\sum\limits_{n_1=0}^2
\Big(\sum\limits_{n_2=0}^2Q(n_1,n_2)\Big)^2.
\end{array}
\]
The sums in the round brackets in \er{fisl3} are the potentials over "horizontal"\, (respectively, "vertical"\,) cycles in the fundamental graph $\cG_*=\Z^2/\G$ (shown in red (respectively, in blue) in Fig.~\ref{fig1}\emph{b}).
\end{remark}

It is convenient to represent the spectral invariants \er{pisl1}, \er{qfZd} in terms of the discrete Fourier transform $\wh Q$ of the potential $Q:\cV_*\to\C$ given by
\[\lb{DFT}
\begin{array}{l}
\displaystyle\wh Q(l)=\frac1{p}\,\sum_{\mn\in\cV_*}e^{-2\pi i\big(\frac{l_1n_1}{p_1}+\ldots+\frac{l_dn_d}{p_d}\big)}Q(\mn),\qq \textrm{where}\qq p=p_1\ldots p_d,\\[14pt]
l=(l_1,\ldots,l_d)\in\cV_*,\qqq \mn=(n_1,\ldots,n_d)\in\cV_*,
\end{array}
\]
$\cV_*$ has the form \er{cV*}.

\begin{corollary}\lb{CSIF}
Let the potential $Q$ be real. Then the periodic spectral invariants \er{pisl1} and the Floquet spectral invariants \er{qfZd}, in terms of the Fourier transform \er{DFT},  have the form
\[\lb{pslF}
\cI_1(Q)=p\,\wh Q(0),\qquad \cI_2(Q)=\frac p2\sum_{\mn\in\cV_*}\big|\wh Q(\mn)\big|^2,
\]
\[\lb{fslF}
\cI^{\mathfrak{e}_j}_{p_j+2}(Q)=\cI_2(Q)+\frac{p\,p_j}2
\sum\limits_{\mn\in\cV_*\atop n_j=0}
\big|\wh Q(\mn)\big|^2,\qq \mn=(n_1,\ldots,n_d), \qq j=1,\ldots,d.
\]
\end{corollary}

\begin{remark}
\emph{i}) The first identity in \er{pslF} holds true for complex potentials $Q$ as well.

\emph{ii}) From \er{qfZd} and \er{fslF}, it follows that
\[\lb{Fpsi}
\sum\limits_{n_1,\ldots,\breve{n}_j,\ldots,n_d}
\Big(\sum\limits_{n_j=0}^{p_j-1}Q(\mn)\Big)^2
\]
or, in terms of the Fourier transform (for real $Q$),
\[\lb{Ffsi}
\sum\limits_{\mn\in\cV_*\atop n_j=0}\big|\wh Q(\mn)\big|^2,\qqq \mn=(n_1,\ldots,n_d), \qqq j=1,\ldots,d,
\]
are also Floquet spectral invariants. The periodic spectral invariants \er{pisl1}, \er{pslF} and the Floquet spectral invariants \er{Fpsi}, \er{Ffsi} for the lattice $\Z^d$ ($d\geq2$) were obtained in \cite{Ka89} using the discrete heat equation for the Schr\"odinger operator $H$ on $\Z^d$ and in \cite{L23,L24} by expanding the dispersion relation for $H$ into the Laurent polynomials. We obtain these invariants using the trace formulas for the Schr\"odinger operators from \cite{KS22}.
\end{remark}

\section{Proof of the main results}
\setcounter{equation}{0}
\lb{Sec3}

\subsection{Trace formulas for Schr\"odinger operators} The proof of our results is based on the following trace formulas for the Schr\"odinger operator (see \cite[Corollary 2.6]{KS22}). For each closed path $\gp=(\be_1,\ldots,\be_n)$ in the modified fundamental graph $\wt\cG_*$, we define the index $\t(\gp)$ and the weight $\o(\gp,Q)$ by
$$
\t(\gp)=\t(\bc),\qqq \o(\gp,Q)=\o(\bc,Q),\qqq \bc=[\be_1,\ldots,\be_n],
$$
where the cycle $\bc$ is an equivalence class of $\gp$, and the index $\t(\bc)$ and the weight $\o(\bc,Q)$ of a cycle $\bc$ are given by \er{cyin} and \er{Wcy}, respectively.

\begin{theorem}\lb{TPG}
Let $H(k)=A(k)+Q$, $k\in\T^d$, be the fiber Schr\"odinger operator defined by \er{Hvt'}, \er{fado} on the fundamental graph $\cG_*$. Then for each $n\in\N$, the trace of $H^n(k)-A^n(k)$ has the form
\[\lb{TrHAn}
\Tr\big(H^n(k)-A^n(k)\big)
=\sum_{\mm\in\Z^d}\gt_n^\mm(Q)\cos\lan\mm,k\ran,
\qqq\gt_n^\mm(Q)=\sum_{\gp\in\wt C_n^\mm\sm C_n^\mm}\o(\gp,Q),
\]
where $C_n^\mm$ and $\wt C_n^\mm$ are the sets of all closed paths of length $n$ and with index $\mm$ in the fundamental graph $\cG_*$ and in the modified fundamental graph $\wt\cG_*$ defined by (\ref{mfg}), respectively, and $\o(\gp,Q)$ is the weight of the closed path $\gp$.
\end{theorem}

\begin{remark}
\emph{i}) The formulas \er{TrHAn} are \emph{trace formulas}, where the traces (the sums of eigenvalues) of the $n$-th power of the fiber operators are expressed in terms of the potential~$Q$ and closed paths of length $n$ in $\wt\cG_*$.

\emph{ii}) The sum over $\mm$ in \er{TrHAn} is finite, since $\wt C_n^\mm=\varnothing$ for all $\mm\in\Z^d$ such that $\|\mm\|>n\t_+$, where $\t_+$ is given in \er{cNe0} (see \cite[Proposition 2.1.\emph{i}]{KS22}).

\emph{iii}) Although in \cite{KS22} the authors considered real potentials $Q$, the trace formulas \er{TrHAn} hold true for complex valued periodic potentials as well.
\end{remark}

\subsection{Spectral invariants for Schr\"odinger operators}
We prove Theorem \ref{TSpI} about complete systems of spectral invariants for the Schr\"odin\-ger operator on periodic graphs.

\medskip

\no \textbf{Proof of Theorem \ref{TSpI}.} \emph{i)} Let $H_s=A+Q_s$ be the Schr\"odinger operators with $\G$-periodic potentials $Q_s$ on a $\G$-periodic graph $\cG$, and let $H_s(k)=A(k)+Q_s$, $k\in\T^d$, be the corresponding Floquet operators, $s=1,2$.

Assume that $\s\big(H_1(k)\big)=\s\big(H_2(k)\big)$ for all $k\in\T^d$. Then
$$
\textstyle\Tr H_1^n(k)=\sum\limits_{j=1}^\n\l_j^n(k,Q_1)=\sum\limits_{j=1}^\n\l_j^n(k,Q_2)=\Tr H_2^n(k),\qqq \forall\,(n,k)\in\N\ts\T^d,
$$
where $\l_j(k,Q_s)$, $j\in\N_\n$, are the eigenvalues of $H_s(k)$. This yields
$$
\Tr\big(H_1^n(k)-A^n(k)\big)=\Tr\big(H_2^n(k)-A^n(k)\big),\qqq \forall\,(n,k)\in\N\ts\T^d.
$$
Then, using \er{TrHAn}, we obtain
\[\lb{summ}
\sum_{\mm\in\Z^d}\gt_n^\mm(Q_1)\cos\lan\mm,k\ran=
\sum_{\mm\in\Z^d}\gt_n^\mm(Q_2)\cos\lan\mm,k\ran, \qqq \forall\,(n,k)\in\N\ts\T^d.
\]
Since $\cos\lan\mm,k\ran$, $\mm\in\Z^d$, are linearly independent functions on $\T^d$, from \er{summ} we deduce that
$$
\gt_n^\mm(Q_1)=\gt_n^\mm(Q_2), \qqq \forall\,(n,\mm)\in\N\ts\Z^d.
$$
Thus, $\gt_n^\mm(Q)$, $(n,\mm)\in\N\ts\Z^d$, are Floquet spectral invariants of the Schr\"odinger operator $H=A+Q$ on $\cG$.

Conversely, let
$$
\gt_n^\mm(Q_1)=\gt_n^\mm(Q_2), \qqq \forall\,(n,\mm)\in\N_\n\ts\Z^d.
$$
Then, using \er{TrHAn}, we get
$$
\Tr H_1^n(k)=\Tr H_2^n(k),\qqq \forall\,(n,k)\in\N_\n\ts\T^d,
$$
which yields (see Remark \ref{RTIn}.\emph{ii}) that $\s\big(H_1(k)\big)=\s\big(H_2(k)\big)$ for all $k\in\T^d$. Thus,
$$
\big\{\gt_n^\mm(Q),\;(n,\mm)\in\N_\n\ts\Z^d\big\}
$$
is a complete system of Floquet spectral invariants of $H$.

Each closed path $\gp$ in the modified fundamental graph $\wt\cG_*$ can be expressed in a unique way in the form $\gp=\gp_o^r$ for some $r\in\N$ and some prime closed path $\gp_o$ in $\wt\cG_*$, and
\[\lb{epc}
|\gp|=r|\gp_o|,\qqq \t(\gp)=r\t(\gp_o),\qqq \o(\gp,Q)=\o^r(\gp_o,Q).
\]

Each prime cycle $\bc$ is an equivalence class consisting of $|\bc|$ distinct prime closed paths. These closed paths are a simple circular shift of one another and, consequently, have the same length $|\bc|$, index $\t(\bc)$ and weight $\o(\bc,Q)$. Then, using \er{epc}, we can rewrite the identity for $\gt_n^\mm(Q)$ in \er{TrHAn} in the form
$$
\gt_n^\mm(Q)=\sum_{\gp\in\wt C_n^\mm\sm C_n^\mm}\o(\gp,Q)=
\sum_{r\in\N,\;\bc\in\wt\cP\sm\cP\atop r|\bc|=n,\,r\t(\bc)=\mm}|\bc|\;\o^r(\bc,Q)=\sum_{r\in\N_n,\;\bc\in\wt\cP\sm\cP\atop r|\bc|=n,\,r\t(\bc)=\mm}\frac nr\;\o^r(\bc,Q)=n\,\cI_n^\mm(Q),
$$
where $\cP$ and $\wt\cP$ are the sets of all prime cycles in the fundamental graph $\cG_*$ and in the modified fundamental graph $\wt\cG_*$, respectively, and $\cI_n^\mm(Q)$ is given by \er{cInm}. This yields that $\cI_n^\mm(Q)$, $(n,\mm)\in\N\ts\Z^d$, are also Floquet spectral invariants of $H$, and the system \er{cssi} is complete.

\emph{ii)} For $k=0$ the identities \er{TrHAn} have the form
$$
\begin{array}{l}
\Tr\big(H^n(0)-A^n(0)\big)=\gt_n(Q), \qqq n\in\N,\qqq \textrm{where} \\[8pt] \gt_n(Q):=\sum\limits_{\mm\in\Z^d}\gt_n^\mm(Q)=
\sum\limits_{\gp\in\wt C_n\sm C_n}\o(\gp,Q),
\end{array}
$$
$C_n$ and $\wt C_n$ are the sets of all closed paths of length $n$ in $\cG_*$ and $\wt\cG_*$, respectively. The remaining part of the proof is similar to the proof of the item \emph{i}). \qq $\Box$

\medskip

We prove Proposition \ref{spcN} about properties of the spectral invariants $\cI_n(Q)$ and $\cI_n^\mm(Q)$, $(n,\mm)\in\N\ts\Z^d$.

\medskip

\no \textbf{Proof of Proposition \ref{spcN}.} \emph{i}) The first statement of this item is a direct consequence of the formulas \er{cInm} and \er{cIn} for the spectral invariants $\cI_n^\mm(Q)$ and $\cI_n(Q)$ and the definition \er{Wcy} of the cycle weight $\o(\bc,Q)$, see also Remark \ref{Rcyw}.\emph{i}.

The identity \er{coPF} follows from \er{cInm} and \er{cIn}.

For each cycle $\bc\in\wt\cP\sm\cP$ the reverse cycle $\ol{\bc}$ is also in $\wt\cP\sm\cP$ and $|\ol{\bc}\,|=|\bc|$, $\t(\ol{\bc}\,)=-\t(\bc)$, $\o(\ol{\bc},Q)=\o(\bc,Q)$, see \er{ininc} and \er{inwe}. This and the definition \er{cInm} of the spectral invariant $\cI_n^\mm(Q)$ yield the identity \er{Npme}.

Each cycle $\bc\in\wt\cP\sm\cP$ has at least one loop edge $\be_v$, $v\in\cV_*$, with zero index, see \er{mfg}. Then, using the definition \er{cyin} of the cycle index, we have
$$
\textstyle\|\t(\bc)\|\leq \sum\limits_{\be\in\bc}\|\t(\be)\|\leq(|\bc|-1)\t_+,
$$
where $\t_+$ is defined in \er{cNe0}. Therefore, for any $\bc\in\wt\cP\sm\cP$ and any $r\in\N$ we obtain
$$
r\|\t(\bc)\|\leq (r|\bc|-r)\t_+\leq (r|\bc|-1)\t_+.
$$
Thus, if $\|\mm\|>(n-1)\t_+$, then there are no $\bc\in\wt\cP\sm\cP$ and $r\in\N_n$ such that $r|\bc|=n$ and $r\t(\bc)=\mm$. Then, using \er{cInm},  we obtain \er{cNe0}.

\emph{ii}) Let $\mm\in\Z^d$ be primitive. Then for any cycle $\bc\in\wt\cP\sm\cP$, if $r\t(\bc)=\mm$ for some $r\in\N_n$, then $r=1$ and the identity \er{cInm} for $\cI_n^\mm(Q)$ has the form \er{Inpm}.

\emph{iii}) Let $n$ be prime. Then for any cycle $\bc\in\wt\cP\sm\cP$, if $r|\bc|=n$ for some $r\in\N_n$, then $r=1$ or $r=n$ and the identities \er{cInm}, \er{cIn} for $\cI_n^\mm(Q)$ and $\cI_n(Q)$ may be written in the form
\[\lb{cInm1}
\cI_n^\mm(Q)=\sum_{\bc\in\wt\cP_n^\mm\sm\cP_n^\mm}\o(\bc,Q)+\sum_{\bc\in\wt\cP_1\sm\cP_1\atop n\t(\bc)=\mm}\frac1n\;\o^n(\bc,Q),\qqq \mm\in\Z^d,
\]
\[\lb{cIn1}
\cI_n(Q)=\sum_{\bc\in\wt\cP_n\sm\cP_n}\o(\bc,Q)+
\sum_{\bc\in\wt\cP_1\sm\cP_1}\frac1n\;\o^n(\bc,Q).
\]
We have $\bc\in\wt\cP_1\sm\cP_1$ if and only if $\bc=[\be_v]$, where $\be_v$, $v\in\cV_*$, is the loop added to the fundamental graph, and $\t(\bc)=0$, $\o(\bc,Q)=Q(v)$. Then \er{cInm1} and \er{cIn1} have the form \er{Ipnm} -- \er{cIpn}. \qq $\Box$

\medskip

We present explicit expressions for the first periodic spectral invariants $\cI_n(Q)$ and Floquet spectral invariants $\cI_n^\mm(Q)$, $\mm\in\Z^d$, $n=1,2,3$. For each cycle $\bc=[v_1,\ldots,v_n]$ (given by the sequence of its vertices $v_1,\ldots,v_n\in\cV_*$) in the fundamental graph $\cG_*=(\cV_*,\cA_*)$, we define the potential $Q(\bc)$ over $\bc$ by
\[\lb{pQoc}
Q(\bc)=Q(v_1)+\ldots+Q(v_n),\qqq \bc=[v_1,\ldots,v_n].
\]

\begin{proposition}\lb{TrTO} For $n=1,2,3$ the Floquet spectral invariants $\cI_n^\mm(Q)$, $\mm\in\Z^d$, and the periodic spectral invariants $\cI_n(Q)$ defined by \er{cInm} and \er{cIn}, respectively, have the form
\[\lb{TrH123}
\begin{array}{l|l}
\mm=0 & \; \mm\neq0\\[6pt]
\cI_1^0(Q)=\sum\limits_{v\in\cV_*}Q(v), & \;\cI_1^\mm(Q)=0,\\[12pt]
\cI_2^0(Q)=\frac12\sum\limits_{v\in\cV_*}Q^2(v), &\; \cI_2^\mm(Q)=\sum\limits_{\bc\in\cP_1^\mm}Q(\bc),\\[12pt]
\cI_3^0(Q)=\frac13\sum\limits_{v\in\cV_*}Q^3(v)+\sum\limits_{\bc\in\cP_2^0}Q(\bc),& \;\cI_3^\mm(Q)=\sum\limits_{\bc\in\cP_1^\mm}Q^2(\bc)+
\sum\limits_{\bc\in\cP_2^\mm}Q(\bc),\qq \mm\notin2\Z^d,\\[12pt]
& \;\cI_3^{2\mm}(Q)=\!\sum\limits_{\bc\in\cP_1^{2\mm}}Q^2(\bc)+
\!\sum\limits_{\bc\in\cP_1^{\mm}}Q(\bc)
+\!\sum\limits_{\bc\in\cP_2^{2\mm}}Q(\bc);
\end{array}
\]
\[\lb{TrI123}
\begin{array}{l}
\cI_1(Q)=\sum\limits_{v\in\cV_*}Q(v), \\[12pt]
\cI_2(Q)=\frac12\sum\limits_{v\in\cV_*}Q^2(v)+
\sum\limits_{\bc\in\cP_1}Q(\bc),\\[12pt]
\cI_3(Q)=\frac13\sum\limits_{v\in\cV_*}Q^3(v)+
\sum\limits_{\bc\in\cP_1}\big(Q^2(\bc)+Q(\bc)\big)+\sum\limits_{\bc\in\cP_2}Q(\bc).
\end{array}
\]
Here $\cP_n$ is the set of all prime cycles of length $n$ in the fundamental graph $\cG_*=(\cV_*,\cA_*)$, and $\cP_n^\mm\ss\cP_n$ is its subset of prime cycles with index $\mm$; the potential $Q(\bc)$ over the cycle $\bc$ is defined by \er{pQoc}.

$\bu$ If the periodic graph $\cG$ has no multiple edges, then
\[\lb{pgwm}
\textstyle\cI_3^0(Q)=\frac13\sum\limits_{v\in\cV_*}Q^3(v)+
\sum\limits_{v\in\cV_*}\vk_vQ(v),
\]
where $\vk_v$ is the degree of the vertex $v$.

$\bu$ If the fundamental graph $\cG_*$ has no loops, then
\[\lb{wolp}
\begin{array}{l}
\cI_2^\mm(Q)=0,\hspace{10mm} \cI_3^\mm(Q)=\sum\limits_{\bc\in\cP_2^\mm}Q(\bc),\hspace{10mm}  \mm\neq0; \\[14pt]
\cI_2(Q)=\cI_2^0(Q)=\frac12\sum\limits_{v\in\cV_*}Q^2(v),
\end{array}
\]
\[\lb{wolp3}
\textstyle \cI_3(Q)=\frac13\sum\limits_{v\in\cV_*}Q^3(v)+\sum\limits_{v\in\cV_*}m_vQ(v),
\hspace{10mm} m_v=\sum\limits_{u\sim v}m_{v,u}^2,
\]
where the sum $\sum\limits_{u\sim v}$ is taken over all vertices $u\in\cV_*$ adjacent to $v$, and $m_{v,u}$ is the multiplicity of the edge $(v,u)\in\cA_*$ in $\cG_*$.

$\bu$ If the fundamental graph $\cG_*$ has no loops and multiple edges, then
\[\lb{wolp1}
\textstyle\cI_3^\mm(Q)=0,\hspace{5mm}  \mm\neq0; \hspace{10mm}
\cI_3(Q)=\cI_3^0(Q)=\frac13\sum\limits_{v\in\cV_*}Q^3(v)+\sum\limits_{v\in\cV_*}\vk_vQ(v).
\]
\end{proposition}

\begin{remark}
Recall that we assume that there are no loops in the periodic graph $\cG$ (see Remark~\ref{Re12}.\emph{ii}). In this case the fundamental  graph $\cG_*$ has no loops with zero index. Note that $\cG_*$ may have loops with non-zero indices.
\end{remark}

\no \textbf{Proof of Proposition \ref{TrTO}.} If $r|\bc|=1$ for some cycle $\bc$ and $r\in\N$, then $r=1$ and $|\bc|=1$. Thus, the identities \er{cInm} for $n=1$ may be written in the form
\[\lb{cInm11}
\textstyle\cI_1^\mm(Q)=\sum\limits_{\bc\in\wt\cP_1^\mm\sm\cP_1^\mm}\o(\bc,Q),\qqq \mm\in\Z^d.
\]
Here and below $\wt\cP_n$ is the set of all prime cycles of length $n$ in the modified fundamental graph $\wt\cG_*$, and $\wt\cP_n^\mm\ss\wt\cP_n$ is its subset of prime cycles with index $\mm$; $\o(\bc,Q)$ is the weight of a cycle $\bc$ defined by \er{Wcy}.

Each cycle $\bc\in\wt\cP_1\sm\cP_1$ has the form $\bc=[\be_v]$, $v\in\cV_*$, where $\be_v$ is the loop added to the fundamental graph, and
$$
\t(\bc)=\t(\be_v)=0,\qqq \o(\bc,Q)=\o(\be_v)=Q(v).
$$
Then the identities \er{cInm11} have the form
$$
\textstyle\cI_1^0(Q)=\sum\limits_{v\in\cV_*}Q(v),\qqq
\cI_1^\mm(Q)=0,\qq \mm\neq0.
$$

Each prime cycle $\bc\in\wt\cP_2\sm\cP_2$ has the form $\bc=[\be_v,\bc_v]$, where $\bc_v\in\cA_*$ is a loop at some vertex $v$ in $\cG_*$ (with weight 1). Due to the definitions \er{cyin} and \er{Wcy} of the cycle index and the cycle weight, we have
$$
\t(\bc)=\t(\be_v)+\t(\bc_v)=\t(\bc_v),\qqq \o(\bc,Q)=\o(\be_v)\o(\bc_v)=Q(v).
$$
We note that $\t(\bc_v)\neq0$, since there are no loops with zero index in the fundamental graph $\cG_*$. Then, using the identities \er{Ipnm}, \er{Ipn0} for $n=2$, we obtain
$$
\textstyle\cI_2^\mm(Q)=\sum\limits_{\bc\in\wt\cP_2^\mm\sm\cP_2^\mm}\o(\bc,Q)=
\sum\limits_{[\bc_v]\in\cP_1^\mm}Q(v),\qqq 0\neq\mm\in\Z^d,
$$
$$
\textstyle\cI_2^0(Q)=\sum\limits_{\bc\in\wt\cP_2^0\sm\cP_2^0}\o(\bc,Q)+\frac12
\sum\limits_{v\in\cV_*}Q^2(v)=\frac12
\sum\limits_{v\in\cV_*}Q^2(v).
$$

Finally, the identities \er{Ipnm}, \er{Ipn0} for $n=3$, have the form
\[\lb{Ipnm3}
\textstyle\cI_3^\mm(Q)=\sum\limits_{\bc\in\wt\cP_3^\mm\sm\cP_3^\mm}\o(\bc,Q),\qqq 0\neq\mm\in\Z^d,
\]
\[\label{Ipn03}
\textstyle\cI_3^0(Q)=\sum\limits_{\bc\in\wt\cP_3^0\sm\cP_3^0}\o(\bc,Q)+
\frac13\sum\limits_{v\in\cV_*}Q^3(v).
\]
Each prime cycle in $\wt\cP_3\sm\cP_3$ has one of the following forms

\no $\bu$ $\bc_1=[\be_v,\be_v,\bc_v]$ or $\bc_2=[\be_v,\bc_v,\bc_v]$, where $\bc_v\in\cA_*$ is a loop at some vertex $v$ in $\cG_*$;

\no $\bu$ $\bc_{3,1}=[\be_1,\be_2,\be_{v_1}]$ or $\bc_{3,2}=[\be_1,\be_{v_2},\be_2]$, where $\be_1,\be_2\in\cA_*$, $\be_1=(v_1,v_2)$, $\be_2=(v_2,v_1)$ for some $v_1,v_2\in\cV_*$ (not necessarily distinct) and $\be_1\neq\be_2$, i.e., $\bc_o:=[\be_1,\be_2]$ is a prime cycle in $\cG_*$ (with weight 1).

\no For these cycles we have
$$
\begin{array}{ll}
\t(\bc_1)=2\t(\be_v)+\t(\bc_v)=\t(\bc_v)\neq0,\qq & \o(\bc_1,Q)=\o(\be_v)\o(\be_v)\o(\bc_v)=Q^2(v);\\[4pt]
\t(\bc_2)=\t(\be_v)+2\t(\bc_v)=2\t(\bc_v)\neq0,\qq & \o(\bc_2,Q)=\o(\be_v)\o(\bc_v)\o(\bc_v)=Q(v);\\[4pt]
\t(\bc_{3,j})=\t(\be_{v_j})+\t(\bc_o)=\t(\bc_o),\; & \o(\bc_{3,j},Q)=\o(\be_{v_j})\o(\bc_o,Q)=Q(v_j),\qq j=1,2.
\end{array}
$$
This yields that \\
$\bu$ if $\mm=0$, then
\[\lb{mmze}
\textstyle\sum\limits_{\bc\in\wt\cP_3^0\sm\cP_3^0}\o(\bc,Q)=
\sum\limits_{\bc_o=[v_1,v_2]\in\cP_2^0}\big(Q(v_1)+Q(v_2)\big)
=\sum\limits_{\bc_o\in\cP_2^0}Q(\bc_o);
\]
$\bu$ if $0\neq\mm\in2\Z^d$, then
\[\lb{mmev}
\textstyle\sum\limits_{\bc\in\wt\cP_3^\mm\sm\cP_3^\mm}\o(\bc,Q)=
\sum\limits_{[\bc_v]\in\cP_1^\mm}Q^2(v)+\sum\limits_{[\bc_v]\in\cP_1^{\mm/2}}Q(v)+ \sum\limits_{\bc_o\in\cP_2^\mm}Q(\bc_0);
\]
$\bu$ if $\mm\notin2\Z^d$, then
\[\lb{mmod}
\textstyle\sum\limits_{\bc\in\wt\cP_3^\mm\sm\cP_3^\mm}\o(\bc,Q)=
\sum\limits_{[\bc_v]\in\cP_1^\mm}Q^2(v)+ \sum\limits_{\bc_o\in\cP_2^\mm}Q(\bc_0).
\]
Substituting \er{mmze} -- \er{mmod} into \er{Ipnm3}, \er{Ipn03}, we obtain
the identities for the Floquet spectral invariants $\cI_3^\mm(Q)$, $\mm\in\Z^d$, given in \er{TrH123}.

The identities \er{TrI123} for the periodic spectral invariants $\cI_n(Q)$, $n=1,2,3$, follow from \er{coPF} and \er{TrH123} and the fact that $\cP_1^0=\varnothing$.

If the periodic graph $\cG$ has no multiple edges, then all cycles of length 2 with zero index in the fundamental graph $\cG_*=(\cV_*,\cA_*)$ are the cycles $\bc_{\be}=[\be,\ol\be\,]$, $\be\in\cA_*$, and $\bc_{\ol\be}=\bc_{\be}$. Then
$$
\textstyle\sum\limits_{\bc\in\cP_2^0}Q(\bc)=
\frac12\sum\limits_{\be=(u,v)\in\cA_*}\big(Q(u)+Q(v)\big)=
\sum\limits_{v\in\cV_*}\vk_vQ(v),
$$
and the invariant $\cI_3^0(Q)$ in \er{TrH123} has the form \er{pgwm}.

If the fundamental graph $\cG_*$ has no loops, then $\cP_1=\cP_1^\mm=\varnothing$ for each $\mm\in\Z^d$ and the identities \er{TrH123}, \er{TrI123} yield \er{wolp} and
\[\lb{I3P2}
\textstyle \cI_3(Q)=\frac13\sum\limits_{v\in\cV_*}Q^3(v)+\sum\limits_{\bc\in\cP_2}Q(\bc).
\]
Let $u,v\in\cV_*$, $u\neq v$, and let $\be_i=(v,u)\in\cA_*$, $i=1,\ldots,m_{v,u}$, where $m_{v,u}$ is the multiplicity of the edge $(v,u)$. Then each cycle of length 2 passing through the vertices $v,u$ has the form $[\be_i,\ol\be_j]$ for some (not necessarily distinct) $i,j=1,\ldots,m_{v,u}$. Thus, the number of such cycles is $m_{v,u}^2$, and
$$
\textstyle \sum\limits_{\bc\in\cP_2}Q(\bc)=\sum\limits_{v\in\cV_*}Q(v)\sum\limits_{u\sim v}m_{v,u}^2.
$$
Combining this with \er{I3P2}, we get \er{wolp3}.

Finally, let the fundamental graph $\cG_*$ have no loops and multiple edges (in this case the periodic graph $\cG$ also has no multiple edges). Then
$$
\textstyle \cP_2^\mm=\varnothing, \;\; \forall\,\mm\in\Z^d\sm\{0\};\qq m_{v,u}=1,\;\; \forall\,(v,u)\in\cA_*;\qq m_v:=\sum\limits_{u\sim v}m_{v,u}^2=\vk_v,\;\; \forall\,v\in\cV_*,
$$
and using \er{pgwm} -- \er{wolp3}, we obtain \er{wolp1}. \qq $\Box$

\medskip

\no \textbf{Proof of Corollary \ref{CSpI}.} The explicit expressions for the invariants $\cI_n(Q)$, $n=1,2,3$, defined by \er{cIn} have already been proved in Proposition \ref{TrTO}. One just needs to put $\vk_v=\vk_o$, $v\in\cV_*$, in the last identity in \er{wolp1}. \qq $\Box$

\medskip

Now we prove Theorem \ref{TLQI} about the linear and quadratic Floquet spectral invariants of the Schr\"odinger operator.

\medskip

\no \textbf{Proof of Theorem \ref{TLQI}.} Let $\mm\in\Z^d$ be primitive, and let $n:=n(\mm)$ be the length of the shortest cycle with index $\mm$ in the fundamental graph $\cG_*=(\cV_*,\cA_*)$. Then, due to \er{Inpm}, the Floquet spectral invariants $\cI_{n+s}^\mm(Q)$, $s=1,2$, have the form
\[\lb{denn}
\cI_{n+s}^\mm(Q)=\sum_{\bc\in\wt\cP_{n+s}^\mm\sm\cP_{n+s}^\mm}\o(\bc,Q),\qqq
s=1,2.
\]
Each cycle in $\wt\cP_{n+s}^\mm\sm\cP_{n+s}^\mm$, $s=1,2$, contains at least one loop $\be_v$, $v\in\cV_*$, added to the fundamental graph $\cG_*$, see \er{mfg}. Recall that each such loop $\be_v$ has the index $\t(\be_v)=0$ and the weight $\o(\be_v)=Q(v)$.

Since the shortest cycle with index $\mm$ in $\cG_*$ has the length $n$, each prime cycle $\bc$ in $\wt\cP_{n+1}^\mm\sm\cP_{n+1}^\mm$ is a union of a cycle $\bc_o=[v_1,\ldots,v_n]\in\cP_n^\mm$ (which is also prime, since its index $\mm$ is primitive) and one of the added loops $\be_{v_j}$, $j\in\N_n$, i.e., $\bc=\bc_o\cup\be_{v_j}$. Then, using the definition \er{Wcy} of the cycle weight $\o(\bc,Q)$, we have
$$
\o(\bc,Q)=\o(\bc_o,Q)\,\o(\be_{v_j})=Q(v_j),
$$
since $\o(\bc_o,Q)=1$ for any cycle $\bc_o$ in $\cG_*$. Therefore, the identity \er{denn} when $s=1$ can be written as
$$
\cI_{n+1}^\mm(Q)=\sum_{\bc_o=[v_1,\ldots,v_n]\in\cP_n^\mm}\sum_{j=1}^nQ(v_j)=\hspace{-3mm}
\sum_{\bc_o=[v_1,\ldots,v_n]\in\cP_n^\mm}\hspace{-3mm}h_1(q_1,\ldots,q_n),\qqq q_j=Q(v_j).
$$
Thus, \er{lfin} is proved.

Similarly, each cycle $\bc$ in $\wt\cP_{n+2}^\mm\sm\cP_{n+2}^\mm$ is \\
\no $\bu$ either a union of a prime cycle $\bc_o=[v_1,\ldots,v_n]\in\cP_n^\mm$ and two added loops $\be_{v_j},\be_{v_l}$, $j,l\in\N_n$ (not necessarily distinct), i.e., $\bc=\bc_o\cup\be_{v_j}\cup\be_{v_l}$, and
$$
\o(\bc,Q)=\o(\bc_o,Q)\,\o(\be_{v_j})\,\o(\be_{v_l})=Q(v_j)Q(v_l);
$$
\no $\bu$ or a union of a prime cycle $\wt\bc_o=[v_1,\ldots,v_{n+1}]\in\cP_{n+1}^\mm$ (if such a cycle exists) and one of the added loops $\be_{v_j}$, $j\in\N_{n+1}$, i.e., $\bc=\wt\bc_o\cup\be_{v_j}$, and $\o(\bc,Q)=\o(\wt\bc_o,Q)\,\o(\be_{v_j})=Q(v_j)$.\\
Then for $s=2$, the identity \er{denn} has the form \er{qfin}.

Let the periodic graph $\cG$ be bipartite. Then there are no odd-length cycles in $\cG$ and, consequently, there are no odd-length cycles with zero index in the fundamental graph $\cG_*$, since each cycle with zero index in $\cG_*$ is obtained by factorization of a cycle with the same length in $\cG$, see Remark \ref{Rein}.\emph{i}. Note that the fundamental graph $\cG_*$ may have odd-length cycles with non-zero indices (see Fig.~\ref{fig1}\emph{b}). We will show that $\cP_{n+1}^\mm=\varnothing$. The proof is by contradiction. We suppose that there exists $\bc_o\in\cP_{n+1}^\mm$. We take a cycle
$\bc\in\cP_n^\mm$ (by the theorem's condition such a cycle $\bc$ exists). Let $v$ be a vertex of the cycle $\bc$, and $v_o$ be a vertex of $\bc_o$. Since the fundamental graph $\cG_*$ is connected, there exists a path $\bp$ from $v$ to $v_o$ in $\cG_*$. We consider the cycle $\wt\bc=[\bc,\bp,\ol\bc_o,\ol\bp\,]$, where $\ol\bc$ denotes the reverse of $\bc$. The length $|\wt\bc|$ and the index $\t(\wt\bc\,)$ of $\wt\bc$ satisfy
$$
\begin{array}{l}
|\wt\bc|=|\bc|+|\bp|+|\ol\bc_o|+|\ol\bp|=n+2|\bp|+n+1=2(n+|\bp|)+1,\\[4pt]
\t(\wt\bc)=\t(\bc)+\t(\bp)+\t(\ol\bc_o)+\t(\ol\bp)=0,
\end{array}
$$
i.e., $\wt\bc$ is an odd-length cycle with zero index in the fundamental graph $\cG_*$. We get a contradiction. Thus, $\cP_{n+1}^\mm=\varnothing$, and the identity \er{qfin} has the form \er{qbfi}. \qq $\Box$

\subsection{Zero and degree potentials} We prove Proposition \ref{Pzpo} about isospectrality to the zero and degree potentials.

\medskip

\no\textbf{Proof of Proposition \ref{Pzpo}.} \emph{i}) If the fundamental graph $\cG_*=(\cV_*,\cA_*)$ has no loops, then, due to Corollary \ref{CSpI},
$$
\textstyle\cI_2(Q)=\frac12\sum\limits_{v\in\cV_*}Q^2(v)
$$
is a periodic spectral invariant of $H$. Since $Q$ is isospectral to the zero potential, then
$$
\textstyle\sum\limits_{v\in\cV_*}Q^2(v)=0.
$$
For real $Q$ this yields $Q\equiv0$.

\emph{ii}) According to Corollary \ref{CSpI},
$$
\textstyle\cI_1(Q)=\sum\limits_{v\in\cV_*}Q(v)
$$
is a periodic spectral invariant of $H$. Since $Q$ is isospectral to the  potential $-\vk$, then
\[\lb{Qisk}
\textstyle\sum\limits_{v\in\cV_*}Q(v)=-\sum\limits_{v\in\cV_*}\vk_v.
\]
Let $f\in\ell^2(\cV_*)\cong\C^\n$, $\n=\#\cV_*$. Recall that $Q$ is real. The quadratic form $\lan H(0)f,f\ran$ is non-positive, since
$$
\s\big(H(0)\big)=\s\big(A(0)+Q\big)=\s\big(A(0)-\vk\big)\qqq\textrm{and}\qqq -\D(0):=A(0)-\vk\leq0.
$$
Here $-\D(0)$ and $H(0)$ are the minus Laplacian and Schr\"odinger operator on the fundamental graph $\cG_*$. For $\1=(1,\ldots,1)\in\C^\n$, using \er{Qisk}, we have
\begin{multline*}
\lan H(0)\1,\1\ran=\lan(A(0)-\vk+\vk+Q)\1,\1\ran\\\textstyle=\lan -\D(0)\1,\1\ran+\sum\limits_{v\in\cV_*}\big(\vk_v+Q(v)\big)=\lan-\D(0)\1,\1\ran=0,
\end{multline*}
i.e., the function $\1$ maximizes the quadratic form $\lan H(0)f,f\ran$. Hence, by the Rayleigh-Ritz theorem (see, e.g., p.176 in \cite{HJ85}), 0 is the largest eigenvalue of $H(0)$ with the eigenfunction $\1$. Thus, $\big(A(0)+Q\big)\1=0$ and consequently $Q\equiv-\vk$.  \qq $\Box$

\section{Examples}
\setcounter{equation}{0}
\lb{Sec4}
We prove Examples \ref{E1DL}, \ref{Enzp}, \ref{EKaL}, \ref{EsqL}, and Corollary \ref{CSIF} about spectral invariants for the Schr\"odinger operator on some concrete periodic graphs.

\medskip

\no \textbf{Proof of Example \ref{E1DL}.} The fundamental graph $\cG_*=\Z/\n\Z$ of the one-dimensional lattice $\Z$ is the cycle of length $\n$ and with index 1, see Fig.~\ref{fig3}\emph{b}. The modified fundamental graph $\wt\cG_*$ is obtained from $\cG_*$ by adding a loop with zero index at each vertex of $\cG_*$ (Fig.~\ref{fig3}\emph{c}). The graph $\wt\cG_*$ has no cycles of length less than or equal to $\n$ with non-zero indices, except the fundamental graph $\cG_*$ itself. Then, due to \er{cInm}, the first identities in \er{Inm1d} follow. These identities with the property \er{coPF} yield the second identities in \er{Inm1d}. \qq $\Box$

\medskip

\no \textbf{Proof of Example \ref{Enzp}.} The fundamental graph $\cG_*=\cG/\Z$ consists of two vertices $v_1,v_2$ and two edges $\be_1,\be_2$ with indices $\t(\be_1)=1$, $\t(\be_2)=0$, see Fig.~\ref{slex1}\emph{b} (and their inverse edges). Thus, $\cG_*$ has only two (prime) cycles of length one (i.e., loops) $\be_1$ and $\ol\be_1$, with cycle indices $\t(\be_1)=1$ and $\t(\ol\be_1)=-1$, respectively. Then, using \er{TrH123}, we obtain
$$
\cI_1^\mm(Q)=0\qq \textrm{for}\qq \mm\neq0, \qqq \textrm{and}\qqq \cI_2^\mm(Q)=0\qq \textrm{for}\qq \mm\neq0,\pm1,
$$
and the Floquet spectral invariants \er{1fsi}. Since the number of the fundamental graph vertices $\n=2$, this and Theorem \ref{TSpI}.\emph{i} (see also Remark \ref{Re27}.\emph{ii}) yield that the invariants \er{1fsi} form a complete system of Floquet spectral invariants of $H$.

Similarly, using \er{TrI123} and Theorem \ref{TSpI}.\emph{ii}, we get that \er{1psi} is a complete system of periodic spectral invariants of $H$, and if $P=(p_1,p_2)$ is isospectral to $Q=(q_1,q_2)$, then
$$
\cI_1(Q)=\cI_1(P),\qqq \cI_2(Q)=\cI_2(P).
$$
This and \er{1psi} yield
$$
\textstyle q_1+q_2=p_1+p_2,\qqq \frac12(q_1^2+q_2^2)+2q_1=\frac12(p_1^2+p_2^2)+2p_1.
$$
This system of equations has two solutions $P_1=Q$ and $P_2=(q_2-2,q_1+2)$. If $Q=0$, then $P_1=0$ and $P_2=(-2,2)$. If $Q=-\vk=-(3,1)$, then $P_1=P_2=-\vk$. \qq $\Box$

\medskip

\no \textbf{Proof of Example \ref{EKaL}.} The fundamental graph $\cG_*$ of the Kagome lattice $\cG$ consists of three vertices $v_1,v_2,v_3$, all of degree 4, and six edges
$$
\be_1=(v_1,v_2),\qq \be_2=(v_2,v_1),\qq \be_3=(v_1,v_3),\qq
\be_4=(v_3,v_1),\qq \be_5=(v_3,v_2),\qq \be_6=(v_2,v_3)
$$
(see Fig.~\ref{fig2}\emph{b}) with indices
$$
\begin{array}{lll}
\t(\be_1)=(0,0), \qqq & \t(\be_3)=(0,0), \qqq & \t(\be_5)=(0,0),\\[2pt]
\t(\be_2)=(1,0), & \t(\be_4)=(0,1), & \t(\be_6)=(1,-1).
\end{array}
$$

Since the Kagome lattice $\cG$ has no multiple edges, and there are no loops in $\cG_*$ and each edge in $\cG_*$ has the multiplicity 2, using \er{TrH123} -- \er{wolp3}, we obtain
\[\lb{siKL}
\begin{array}{ll}
\cI_1(Q)=\cI_1^0(Q)=\sum\limits_{s=1}^3q_s, \qqq & \cI_2(Q)=\cI_2^0(Q)=\frac12\sum\limits_{s=1}^3q_s^2,\\[12pt]
\cI_1^\mm(Q)=0, \qq\qq  \cI_2^\mm(Q)=0,&
\cI_3^\mm(Q)=\sum\limits_{\bc\in\cP_2^\mm}Q(\bc),\qqq \mm\neq0,\\[12pt]
\cI_3^0(Q)=\frac13\sum\limits_{s=1}^3q_s^3+4\sum\limits_{s=1}^3q_s, \qqq &
\cI_3(Q)=\frac13\sum\limits_{s=1}^3q_s^3+8\sum\limits_{s=1}^3q_s,
\end{array}
\]
where $q_s=Q(v_s)$, $s=1,2,3$; $\cP_2^\mm$ is the set of all prime cycles of length 2 and with index $\mm$ in $\cG_*$, and $Q(\bc)$ is defined by \er{pQoc}.

All cycles of length 2 and with non-zero indices in $\cG_*$ are the prime cycles
$$
\bc_1=[\be_1,\be_2], \qqq \bc_2=[\be_3,\be_4], \qqq \bc_3=[\be_5,\be_6]
$$
with indices
$$
\t(\bc_1)=(1,0),\qqq \t(\bc_2)=(0,1),\qqq \t(\bc_3)=(1,-1),
$$
and their reverse cycles $\ol\bc_j$ with indices $\t(\ol\bc_j)=-\t(\bc_j)$, $j=1,2,3$. Then, using the identity for $\cI_3^{\mm}(Q)$ in \er{siKL}, we obtain
$$
\cI_3^{\mm}(Q)=0, \qqq \textrm{if}\qqq \mm\neq0,(\pm1,0),(0,\pm1),(\pm1,\mp1),
$$
and the Floquet spectral invariants \er{fiKl}.

Due to Theorem \ref{TSpI}.\emph{ii}, the invariants $\cI_n(Q)$, $n=1,2,3$, form a complete system of the periodic spectral invariants of $H$.

The invariants \er{fiKl} form a complete system of the Floquet spectral  invariants of $H$, since they determine the potential $Q$ uniquely. \qq $\Box$

\medskip

\no \textbf{Proof of Example \ref{EsqL}.} \emph{i}) The fundamental graph $\cG_*=\Z^d/\G=(\cV_*,\cA_*)$ of the lattice $\Z^d$ is a regular graph of degree $2d$ without loops (for $d=2$ and $p_1=p_2=3$ see Fig.~\ref{fig1}\emph{b}). Then, due to Corollary~\ref{CSpI}, the periodic spectral invariants $\cI_n(Q)$, $n=1,2$, defined by \er{cIn} have the form \er{pisl1}. Moreover, if $p_j\geq3$ for each $j\in\N_d$ (i.e., $d_o=0$), then $\cG_*$ has no multiple edges, and using Corollary~\ref{CSpI}, we obtain
$$
\textstyle \cI_3(Q)=\frac13\sum\limits_{\mn\in\cV_*}Q^3(\mn)+2 d\sum\limits_{\mn\in\cV_*}Q(\mn)=\frac13\sum\limits_{\mn\in\cV_*}Q^3(\mn)+
2d\,\cI_1(Q).
$$
If $p_j=2$ for some $j\in\N_d$, then $\cG_*$ has multiple edges (see Fig.~\ref{fig1}\emph{c}). Let $d_o\leq d$ be the number of such $p_1,\ldots,p_d$. Then for each $\mn\in\cV_*$ there are $d_o$ vertices $l\in\cV_*$ such that the edge $(\mn,l)\in\cA_*$ has the multiplicity $m_{\mn,l}=2$ and $2(d-d_o)$ vertices $l'$ such that the edge $(\mn,l')\in\cA_*$ has the multiplicity $m_{\mn,l'}=1$. Thus, we have
$$
\textstyle m_\mn:=\sum\limits_{l\sim\mn}m_{\mn,l}^2=d_o\cdot2^2+2(d-d_o)=2(d+d_o),
$$
where the sum is taken over all vertices $l\in\cV_*$ adjacent to $\mn$,
and, using \er{wolp3}, we obtain
\begin{multline*}
\textstyle \cI_3(Q)=\frac13\sum\limits_{\mn\in\cV_*}Q^3(\mn)+\sum\limits_{\mn\in\cV_*}m_\mn Q(\mn)\\
\textstyle =\frac13\sum\limits_{\mn\in\cV_*}Q^3(\mn)+\sum\limits_{\mn\in\cV_*}2\,(d+d_o)Q(\mn)=
\frac13\sum\limits_{\mn\in\cV_*}Q^3(\mn)+2\,(d+d_o)\cI_1(Q).
\end{multline*}
Thus, \er{pisl2} is proved.

\emph{ii}) The length of the shortest cycle with index $\mathfrak{e}_j$ in the fundamental graph $\cG_*$ is $p_j$, $j\in\N_d$. We derive the invariants \er{lfZd}, \er{qfZd} for $j=1$. The others can be obtained similarly. Recall that we identify the vertices of $\cG_*$ with the corresponding vertices of $\Z^d$ in the fundamental cell $\Omega$ of the period lattice $\G$, see Fig.~\ref{fig1}. All (prime) cycles of length $p_1$ and with index $\mathfrak{e}_1$ in $\cG_*$ are
$$
\begin{array}{ll}
\big[(0,n),(1,n),\ldots,(p_1-1,n)\big],\qqq & n\in\cN:=\cN_{p_2}\ts\ldots\ts\cN_{p_d},\\[6pt]
\cN_{p_l}=\{0,1,\ldots,p_l-1\},\qqq & l=2,\ldots,d,
\end{array}
$$
(given by the sequences of their vertices). The lattice $\Z^d$ is a bipartite graph. Then, using \er{lfin}, \er{qbfi}, we obtain the Floquet spectral invariants
\[\lb{lqfs}
\textstyle\cI^{\mathfrak{e}_1}_{p_1+s}(Q)=\sum\limits_{n\in\cN}h_s\big(Q(0,n),Q(1,n),\ldots,Q(p_1-1,n)\big),\qqq s=1,2,
\]
where $h_s(x_1,\ldots,x_p)$ is the complete homogeneous symmetric polynomial of degree $s$ in variables $x_1,\ldots,x_p$. For $s=1$, the invariant \er{lqfs} has the form
$$
\textstyle\cI^{\mathfrak{e}_1}_{p_1+1}(Q)=
\sum\limits_{\mn\in\cV_*}Q(\mn)=\cI_1(Q).
$$

Using the Newton's identity $2h_2=h_1^2+\phi_2$, where $\phi_2(x_1,\ldots,x_p)=x_1^2+\ldots+x_p^2$ is the second power sum symmetric polynomial in variables $x_1,\ldots,x_p$, we write the invariant \er{lqfs} for $s=2$ in the form
\begin{multline*}
\textstyle2\,\cI^{\mathfrak{e}_1}_{p_1+2}(Q)=
\sum\limits_{n\in\cN}h_1^2\big(Q(0,n),\ldots,Q(p_1-1,n)\big)+
\sum\limits_{n\in\cN}\phi_2\big(Q(0,n),\ldots,Q(p_1-1,n)\big)\\
\textstyle=\sum\limits_{n\in\cN}h_1^2\big(Q(0,n),\ldots,Q(p_1-1,n)\big)+
\sum\limits_{\mn\in\cV_*}Q^2(\mn)\\
\textstyle=\sum\limits_{n\in\cN}h_1^2\big(Q(0,n),\ldots,Q(p_1-1,n)\big)+
2\cI_2(Q),
\end{multline*}
which yields \er{qfZd} when $j=1$. \qq $\Box$

\medskip

\no \textbf{Proof of Corollary \ref{CSIF}.}  Using the definition \er{DFT} of the discrete Fourier transform and Parseval's identity
$$
\textstyle\sum\limits_{\mn\in\cV_*}\big|Q(\mn)\big|^2=p\sum\limits_{\mn\in\cV_*}\big|\wh Q(\mn)\big|^2,
$$
we deduce that the invariants \er{pisl1} have the form \er{pslF}.

Now, we obtain the invariant \er{fslF} for $j=1$. The other ones can be derived similarly. Using \er{DFT} and the identity
$$
\textstyle\sum\limits_{j=0}^{N-1}e^{\frac{2\pi i}N\,(n-n')\,j}=N\,\d_{n,n'},\qqq n,n'=0,1,\ldots,N-1,
$$
where $\d_{n,n'}$ is the Kronecker delta, we obtain
\begin{multline*}
p_1\sum\limits_{l=(0,l_2,\ldots,l_d)\in\cV_*}\big|\wh Q(l)\big|^2=p_1
\sum\limits_{l=(0,l_2,\ldots,l_d)\in\cV_*}\wh Q(l)\,\ol{\wh Q(l)}\\=\frac{p_1}{p^2}\sum\limits_{(0,l_2,\ldots,l_d)\in\cV_*}
\sum_{\mn=(n_1,n_2,\ldots,n_d)\in\cV_*\atop\mn'=(n'_1,n'_2,\ldots,n'_d)\in\cV_*}
e^{-2\pi i\big(\frac{l_2(n_2-n'_2)}{p_2}+\ldots+\frac{l_d(n_d-n'_d)}{p_d}\big)}
Q(\mn)\ol{Q(\mn')}\\
=\frac1p\sum_{(n_1,n_2,\ldots,n_d)\in\cV_*}\sum_{n_1'=0}^{p_1-1}
Q(n_1,n_2,\ldots,n_d)\ol{Q(n'_1,n_2,\ldots,n_d)}\\=
\frac1p\sum_{n_2,\ldots,n_d}\hspace{-2mm}
\bigg(\sum_{n_1,n_1'=0}^{p_1-1}\hspace{-2mm}
Q(n_1,n_2,\ldots,n_d)\ol{Q(n'_1,n_2,\ldots,n_d)}\bigg)=
\frac1p\sum\limits_{n_2,\ldots,n_d}
\Big|\sum\limits_{n_1=0}^{p_1-1}Q(n_1,n_2,\ldots,n_d)\Big|^2.
\end{multline*}
Then for $j=1$ and real $Q$ the invariant \er{qfZd} has the form \er{fslF}.\qq $\Box$

\medskip

\textbf{Acknowledgments.}  This work is supported by the Russian Science Foundation (project No. 25-21-00157). I would like to thank the referees for thoughtful comments that helped me to improve the manuscript.

\medskip

\textbf{Data availability}  No datasets were generated or analysed during the current study.

\medskip

\textbf{Conflict of interest} None declared.

\end{document}